\documentclass[10pt]{article}

\usepackage{amsmath,amssymb}
\usepackage[a4paper]{geometry}

\newtheorem{Lemme}{Lemma}[section] 
\newtheorem{theorem}[Lemme]{Theorem}
\newtheorem{proposition}[Lemme]{Proposition}
\newtheorem{lemma}[Lemme]{Lemma}

\newtheorem{definition}[Lemme]{Definition}

\def\Box{\leavevmode\vrule height 5pt width 4pt depth 0pt\relax}
\date{\today}

\begin{document}

\begin{center}
{\Large \textbf{Optimal control of non-stationary second grade fluids with Navier-slip boundary conditions}}

\vspace{1.3 cm}
\textsc{ nadir arada}{\footnote{Centro de Matem\'atica e Aplica\c{c}\~oes, Departamento de Matem\'atica, Faculdade de Ci\^encias e Tecnologia, Universidade Nova de Lisboa, Portugal. E-mail: naar{\char'100}fct.unl.pt.},} 
  \hspace{10 mm}
\textsc{Fernanda Cipriano}{\footnote{Centro de Matem\'atica e Aplica\c{c}\~oes, Departamento de Matem\'atica, Faculdade de Ci\^encias e Tecnologia, Universidade Nova de Lisboa, Portugal. E-mail: cipriano{\char'100}fct.unl.pt.}}



\end{center}

\vspace{ 1.3 cm}
\noindent \begin{abstract} \noindent We consider optimal control problems governed by systems describing the unsteady flows of an incompressible second grade fluid with Navier-slip boundary conditions. We prove the existence of
an optimal solution and derive the corresponding necessary optimality conditions.\vspace{2mm}\\
  {\bf Key words.} Optimal control, second grade fluid, Navier-slip boundary conditions, necessary optimality conditions.\vspace{3mm}\\
     {\bf AMS Subject Classification.} $49$K$20$, $76$D$55$, $76$A$05$.\vspace{3mm}
\end{abstract} 
\section{Introduction and main results}
\setcounter{equation}{0}
The paper is devoted to the study  of an optimal control problem associated  with a non-stationary viscous, incompressible, second grade fluid.  The state equation is given by
 	\begin{equation}\label{equation_etat}
	\left\{ \begin{array}{ll}
	\frac{\partial }{\partial t}\left( \boldsymbol y-\alpha\Delta  \boldsymbol y
	\right)
	-\nu \Delta  \boldsymbol y+
	\boldsymbol{curl}\left( \boldsymbol y-\alpha\Delta  \boldsymbol y
	\right)\times  \boldsymbol y+\nabla \pi= \boldsymbol u& \quad\mbox{in} \ Q,\vspace{3mm} \\
             \mathrm{div} \,  \boldsymbol y=0&\quad \mbox{in} \ Q,\vspace{3mm}\\
	\boldsymbol y\cdot \boldsymbol n=0,\qquad \left(\boldsymbol n\cdot D\boldsymbol 
	y\right)\cdot \boldsymbol\tau=0&\quad 
	\mbox{on} \ \Sigma,\vspace{3mm} \\
	\boldsymbol y(0)=\boldsymbol y_0&\quad
	\mbox{in} \ \Omega,	
	\end{array}\right.\end{equation}
where $y$ is the velocity field, $\alpha\geq 0$ is a viscoelastic parameter, $\nu>0$ is the viscosity of the fluid,
$\pi$ is the hydrodynamic pressure,   $u$ is a distributed control, $Q=]0,T[\times \Omega$  where $T$ is a fixed positive number and $\Omega\subset \mathbb{R}^2$ is a bounded domain with boundary $\Gamma$,  $\Sigma=]0,T[\times \Gamma$, $n=(n_1,n_2)$ and $\tau=(-n_2,n_1)$ are the unit normal and tangent vectors, respectively, to the boundary $\Gamma$, $Dy=\tfrac{1}{2}\left(\nabla y+\nabla y^\top\right)$ is the symmetric part of the velocity gradient and $y_0\in H^3(\Omega)$ satisfies the incompressibility condition $(\ref{equation_etat})_{2}$ and the boundary conditions $(\ref{equation_etat})_{3}$.  As this equation is set in dimension two, the vector $\boldsymbol y$ is written in the form $\boldsymbol y=(y\equiv (y_1,y_2),0)$ in order to define the curl and the vector product, $\boldsymbol {curl} \, \boldsymbol y=(0,0, \mathrm{curl} \,  y)$ with $ \mathrm{curl} \,  y=\tfrac{\partial y_2}{\partial x_1}-\tfrac{\partial y_1}{\partial x_2}$.  \vspace{1mm}\\
In the inviscid case ($\nu=0$), the second-grade fluid equations are called $\alpha$-Euler equations. Initially proposed as a regularization of the incompressible Euler equations, they are geometrically significant and have been interpreted as a model of turbulence (cf. \cite{HMR981} and \cite{HMR98}). They also inspired another variant, called the $\alpha$-Navier-Stokes equations that turned out to be very relevant in turbulence modeling (cf. \cite{FHT1}, \cite{FHT2} and the references therein). These equations contain the regularizing term $-\nu\Delta\left(\boldsymbol y-\alpha \Delta\boldsymbol y\right)$ instead of $\nu\Delta\boldsymbol y$, making the dissipation stronger and the problem much easier to solve than in the case of second-grade fluids. When $\alpha=0$, the $\alpha$-Navier-Stokes and the second grade fluid equations are equivalent to the Navier-Stokes equations. \vspace{1mm}\\
System $(\ref{equation_etat})_{1,2,4}$ can be supplemented with different kinds of boundary conditions. The case of Dirichlet boundary conditions have received a lot of attention. It was systematically studied for the first time in \cite{O81} and \cite{CO84} for both steady and unsteady cases. 
A Galerkin's method in the basis of the eigenfunctions of the operator $\boldsymbol {curl}(\boldsymbol {curl}(\boldsymbol y-\alpha\Delta \boldsymbol y))$ was especially designed to decompose the problem into a mixed elliptic-hyperbolic type, looking for the velocity $\boldsymbol y$ as a solution of a Stokes-like system coupled to a transport equation satisfied by $\boldsymbol {curl}\left( \boldsymbol y-\alpha\Delta  \boldsymbol y
	\right)$. This approach is optimal in the sense that allows the authors to fully solve the two dimensional problem for both steady and unsteady cases and to {\it automatically} recover  $H^3$ in space regularity. Much work has been done since these pioneering results and, without ambition for completeness, we cite the extensions in \cite{GS94} and \cite{CG97} where global existence for small initial data in three dimensions was established, the former work using a Schauder fixed point argument while the latter considers the decomposition method on the system of Galerkin equations previously mentionned.\vspace{1mm}\\
The case of second grade fluids with Navier boundary conditions has also been  particularly considered in the literature and was studied in \cite{BR03}. These conditions are known to deeply modify the properties of the equations, generating additional difficulties related with boundary terms to be correctly handled. In return, some mathematical aspects turned out to be more easily treatable. This is for example the case when studying the controllability of the Navier-Stokes equations (see \cite{C96}). This is also the case when dealing with the inviscid limit of their solutions. Indeed, it is well known that the solution to the Navier-Stokes equations with Navier boundary conditions converges, as $\nu$ tends to zero, to a solution to the Euler equations, while no similar conclusion can be reached when dealing with Dirichlet boundary conditions, responsible for the formation of boundary layers 
(cf. \cite{CMR98}, \cite{IP06}, \cite{K06}, \cite{LNP05}). Similar considerations apply 
  when analyzing the asymptotic bahavior of the solutions of second-grade fluid equations when the elastic response $\alpha$ and/or the viscosity $\nu$ vanish (cf. \cite{BILN12}, \cite{LT10}, \cite{LNTZ15}).\vspace{1mm}\\
Our objective here is to match the velocity field to a given target field  $y_d\in L^2(Q)$ and the optimal control problem reads as 
	$$(P) \quad
	 \left\{\begin{array}{ll}\mbox{minimize} & \displaystyle J(u,y)=\tfrac{1}{2}\int_Q\left|y-y_d\right|^2\,dxdt+\tfrac{\lambda}{2}
	\int_Q\left|u\right|^2\,dxdt\vspace{1mm}\\
	\mbox{subject to} & (u,y)\in U_{ad}\times L^\infty\left(0,T;H^3(\Omega)\right) \ \mbox{satisfies} \
 (\ref{equation_etat}) \ \mbox{for some} \ \pi\in L^2(\Omega),\end{array}\right.$$
where  $\lambda\geq 0$  and $U_{ad}$, the set of admissible controls, is a nonempty bounded closed convex subset of 
$L^2(I;H({\rm curl};\Omega))=\left\{v\in L^2(Q)\mid {\rm curl} \, v\in L^2(Q)\right\}$. \vspace{2mm}\\
As is well known when dealing with the optimal control of non-Newtonian fluids, the strong nonlinearity in the state equation induces some additional issues  (see {\it e.g.} \cite{A12}, \cite{A13}, \cite{A14}). The first of these difficulties arises when studying the differentiability of the control-to-state mapping and, consequently, when exploring the solvability of the associated linearized equation. In the two dimensional case, this equation can be effectively addressed if the coefficients in the main part of the linearized operator are regular: by expanding the system in the special Galerkin basis used to study the state equation, we may prove the existence and uniqueness of a regular solution without restraining the data. This method, particularly interesting in this situation, shows its limits if the required regularity property is not available: this is the case if the control variables (as well as the initial data) are not regular,
since the coefficients involve the state variable. The second main difficulty, exacerbated by the complexity of the associated differential operators, is encountered when dealing with the adjoint equation.
For both problems, these drawbacks can be overcome by constructing the solutions using the basis of eigenfunctions of the Stokes operator designed in \cite{CMR98} in the case of Navier boundary conditions. Unlike the case of Dirichlet boundary conditions, we are able to derive some corresponding $H^2$ in space a priori estimates and some related time derivative a priori estimates, and it turns out that theses estimates are sufficient to carry out our analysis and prove our main result:
\begin{theorem} \label{main_1} The optimal control problem $(P)$ admits at least one solution $(\bar u,\bar y)$. Moreover, there exists $\bar p\in L^\infty\left(I;H^2(\Omega)\right)$ with 
$\tfrac{\partial \bar p}{\partial t}\in L^2(I;H^1(\Omega))$, unique solution of the following adjoint equation
\begin{equation}\label{adj_opt_eq_alpha}
	\left\{ \begin{array}{lll}
 -\displaystyle \tfrac{\partial }{\partial t}	
	\left(\bar{\boldsymbol p}-\alpha \Delta \bar{\boldsymbol p}\right)-\nu \Delta  \bar{\boldsymbol p}-
	\boldsymbol{curl}\left(\bar{\boldsymbol y}-\alpha \Delta \bar{\boldsymbol y}\right)\times
	\bar{\boldsymbol p}\vspace{2mm} \\
	\hspace{3.5cm}+\boldsymbol{curl}
	\left(\bar{\boldsymbol y}\times \bar{\boldsymbol p}-\alpha \Delta 
	\left(\bar{\boldsymbol y}\times \bar{\boldsymbol p}\right)\right)+\nabla \pi
	= \bar{\boldsymbol y}-\boldsymbol y_d& \quad\mbox{in} \ Q,\vspace{3mm} \\
             \mathrm{div} \,  \bar{\boldsymbol p}=0& \quad\mbox{in} \ Q,\vspace{2mm}\\
                \bar{\boldsymbol p}\cdot \boldsymbol n=0,\qquad \left(\boldsymbol n\cdot D\bar{\boldsymbol p}\right)\cdot \boldsymbol\tau=0& \quad\mbox{on} \ \Sigma,
                \vspace{2mm}\\
    \bar{\boldsymbol p}(T)=0& \quad\mbox{in} \ \Omega,
                
                \end{array}\right.\end{equation}
and satisfying the optimality condition 
	\begin{equation}\label{opt_control_alpha}\int_Q\left(\bar{ p}+\lambda\bar{ u},v-\bar{ u}\right)dxdt\geq 0 \qquad 
	\mbox{for all} \ v\in U_{ad}.\end{equation}
\end{theorem}
The plan of the present paper is as follows. The main results are stated in Section 1. Notation and preliminary results are given in Section 2. Section 3 is devoted to the solvability of the state equation and to the derivation of some corresponding a priori estimates. In Section 4, we establish  existence and uniqueness results for the linearized state equation and  analyze the Lipschitz continuity and the G\^ateaux differentiability of  the control-to-state mapping. The solvability of the adjoint equation is considered in Section 5 and the main results are proved in Section 6. 
\section{Notation and preliminary results}
\setcounter{equation}{0}
Throughout the paper $\Omega$ is a bounded and simply connected domain in $\mathbb{R}^2$. The boundary of $\Omega$ is denoted by $\Gamma$ and is sufficiently regular. We will denote by $I$ the interval $]0,T[$. For $u, v\in \mathbb{R}^{2}$, we define the scalar product by	$u\cdot v=\sum_{i=1}^2 u_iv_i$.
For $\eta, \zeta\in \mathbb{R}^{2\times 2}$, we define the scalar product by	$\eta:\zeta=\sum_{i,j=1}^2 \eta_{ij}\zeta_{ij}$.
 We will also use the following notation
	$$\left(u,v\right)=\displaystyle\int_\Omega u(x)\cdot v(x)\,dx,\qquad 
	\left(\eta,\zeta\right)=\displaystyle\int_\Omega \eta(x):\zeta(x)\,dx.$$
 The standard Sobolev spaces are denoted by $W^{k,p}(\Omega)$ ($k\in \mathbb{N}_0$ and $1<p<\infty$), and their norms by $\|\cdot\|_{k,p}$ and 
$\|\cdot\|_{0,p}\equiv \|\cdot\|_{p}$. We set $W^{k,2}(\Omega)\equiv H^k(\Omega)$ and $\|\cdot\|_{k,2}\equiv \|\cdot\|_{H^k}$. Similarly, the norm in $L^2(Q)$ will be denoted by $\left\|\cdot\right\|_{2,Q}$. 
Since many of the quantities occuring in the paper are vector-valued functions, the notation will be abridged for the sake of brevity and we will omit the space dimension in the function space notation. (The meaning should be clear from the context.) \vspace{1mm}\\
We also introduce the Hilbert space
	$$H({\rm curl};\Omega)=\left\{v\in L^2(\Omega)\mid {\rm curl} \, v\in L^2(\Omega)\right\}$$
and in order to eliminate the pressure in the weak formulation of the state, the linearized state and  the adjoint equations, we consider the following divergence-free spaces
$$H=\left\{v\in L^2(\Omega)\mid \mathrm{div} \, v=0 \ \text{ in }
\Omega \ \mbox{ and } \ v\cdot n=0 \ \mbox{ on }\Gamma\right\},$$
$$V=\left\{v\in H^1(\Omega)\mid \mathrm{div} \, v=0 \
 \mbox{ in } \ \Omega\mbox{ and } \ v\cdot n=0 \ \text{ on } \ \Gamma
 \right\},$$
$$W=\left\{v\in V\cap H^2(\Omega)\mid \left(n\cdot Dv\right)\cdot \tau
=0 \ \ \mbox{on} \ \Gamma\right\}.$$
In the sequel, we set
	$$\sigma(y)=y-\alpha \Delta y \qquad \mbox{for} \ y\in H^2(\Omega)$$
and denote by $\mathbb{P}:L^2(\Omega)\longrightarrow H$, the 
Helmholtz projector in $L^2(\Omega)$ and set $\mathbb{A}=-\mathbb{P}\Delta$. It is well know that $\mathbb{P}$ is a linear bounded operator and that is characterized by the equality $\mathbb P y=\tilde y$, where $\tilde y$ is given by the Helmholtz decomposition 
	$$y=\tilde y+\nabla \phi, \qquad \tilde y\in H 
	\quad \mbox{and} \quad \phi\in H^1(\Omega).$$
Let us now present some useful results. The first one is fundamental and 
 deals with a boundary identity related with the Navier-slip boundary conditions.  It states in particular that the trace of $\mathrm{curl}\, y$ is a linear function on $y$. (See Proposition 1 in \cite{BR03}.)
\begin{lemma} \label{curl_trace}Let $y\in W$. Then, the following identity holds
	$$\mathrm{curl}\, y\big |_\Gamma=y\cdot g\big|_\Gamma \qquad \mbox{where} \ g=2\tfrac{\partial n}{\partial \tau},$$
with $\tfrac{\partial }{\partial \tau}=n_1\tfrac{\partial }{\partial x_2}-n_2\tfrac{\partial }{\partial x_1}$.\end{lemma}
The next two lemmas will be useful when dealing with a priori estimates for the state, linearized and adjoint state equations. We confer to Lemma 5 in \cite{BI06}, Propositions 3 in \cite{BR03} and Lemma 2.1 in \cite{CG97} where similar results are established.
\begin{lemma} \label{sigma_psigma}Let $y\in W\cap H^3(\Omega)$. Then, the following estimate holds
	\begin{equation}\label{sigma_psigma2}\left\|\Delta y+\mathbb{A}y\right\|_{H^1}\leq 
	c\left\|y\right\|_{H^2},\end{equation}
where $c$ is a positive constant only depending on $\Omega$.
\end{lemma}
{\bf Proof.} From the definition of $\mathbb{A}$, there exists $\phi\in H^1(\Omega)$ such that
	$$\Delta y+\mathbb{A} y=\nabla \phi$$
and thus 
	$$\Delta \phi=\mathrm{div}\left(\nabla \phi\right)=
	\mathrm{div}\left(\Delta y+\mathbb{A}y\right)=0.$$
On the other hand, by taking into account Lemma \ref{curl_trace} we obtain
	$$\begin{array}{ll}
	\tfrac{\partial \phi}{\partial n}\big|_{\Gamma}&=
	n\cdot \nabla \phi\big|_{\Gamma}=
	n\cdot \left(\Delta y+\mathbb{A}y\right)\big|_{\Gamma}\vspace{3mm}\\
	&=n\cdot \Delta y\big|_{\Gamma}=\, n\cdot \Delta y\big|_{\Gamma}\vspace{3mm}\\
	&=\, \boldsymbol n\cdot \Delta \boldsymbol y\big|_{\Gamma}
	=-\, \boldsymbol n\cdot \boldsymbol{curl}\left(\boldsymbol{curl}\,\boldsymbol y\right)\big|_{\Gamma}\vspace{3mm}\\
	&=-\, \tfrac{\partial}{\partial \tau}
	\left(y\cdot g\right)\big|_{\Gamma}.\end{array}$$
 Since $\tfrac{\partial}{\partial \tau}
	\left(y\cdot g\right)$ is well defined, the result follows by using standard trace estimates and the regularity theory for elliptic equations with Neuman boundary conditions.$\hfill\Box$ 
\begin{lemma} Let $y\in W\cap H^3(\Omega)$. Then, the following estimates hold
\begin{equation}\label{yh_sigma}\left\|y\right\|_{H^2}\leq 
	c\left(\left\| y\right\|_2+\left\|\mathbb{A} y\right\|_2\right), \end{equation}
where $c$ is a positive constant only depending on $\Omega$.\end{lemma}
{\bf Proof.} Let us first recall that for $f\in H^m(\Omega)$, $m\in \mathbb{N}$, the following problem
	$$\left\{\begin{array}{ll}-\Delta h+h+\nabla \pi=f& \quad \mbox{in} \ \Omega,\vspace{2mm}\\
	\mathrm{div} \, h=0& \quad\mbox{in} \ \Omega,\vspace{2mm}\\
	h\cdot n=0, \qquad (n\cdot Dh)\cdot \tau=0& \quad\mbox{on} \ \Gamma,\vspace{2mm}\\\end{array}\right.$$
admits a unique (up to a constant for $\pi$) solution $(h,\pi)\in H^{m+2}(\Omega)\times H^m(\Omega)$ (see \cite{S73}). Classical arguments show that
	$$\left\|h\right\|_2^2+2 \left\|Dh\right\|_2^2=\left(f,h\right)\leq \left\|f\right\|_2\left\|h\right\|_2$$
yielding 
	\begin{equation}\label{est_g_2_h1}\left\|h\right\|_2\leq \left\|f\right\|_2 \quad \mbox{and} \quad 
	\left\|Dh\right\|_2\leq \tfrac{1}{\sqrt{2}}\left\|f\right\|_2.\end{equation}
On the other hand, due to the regularity results for the Stokes system, we have
	$$\left\|h\right\|_{H^2}\leq 
	c\left\|f-h\right\|_2.$$  
Taking into account (\ref{est_g_2_h1}), we deduce that
	$$\left\|h\right\|_{H^2}\leq c\left(\left\|f\right\|_2+\left\|h\right\|_2\right)\leq c
	\left\|f\right\|_2$$
and the claimed result follows by setting $f=y+\mathbb{A} y$. $\hfill \Box$
\section{State equation}
\setcounter{equation}{0}
In the present section, we state some well known existence and uniqueness results related with the state equation. We first recall an identity  relating the nonlinear term in (\ref{equation_etat}) to the classical trilinear form
	$$b(\phi,z,y)=\left(\phi\cdot \nabla z,y\right)$$  
used in the Euler and Navier-Stokes equations.
For $y,z \in  W\cap H^3(\Omega)$ and $\phi\in  V$, we have
	\begin{align}\label{trilinear_state}\left(\boldsymbol {curl}\,\sigma(\boldsymbol y)\times \boldsymbol z,\boldsymbol\phi\right)&=\left(\boldsymbol {curl}\,\sigma(\boldsymbol y),\mathbf  z \times \boldsymbol\phi\right)\nonumber\\
	&=\displaystyle \left(\sigma(\boldsymbol y),\boldsymbol {curl}\left(\mathbf  z \times \boldsymbol\phi\right)\right)+\int_\Gamma
	\sigma(\boldsymbol y)\times \left(\mathbf  z \times \boldsymbol\phi\right)\cdot \boldsymbol n\,dS\nonumber\\
	&= \left(\sigma(\boldsymbol y),\boldsymbol {curl}\left(\mathbf  z \times \boldsymbol\phi\right)\right)\nonumber\\
&= \left(\sigma(\boldsymbol y),\left( \mathrm{div}\, \boldsymbol\phi\right)\boldsymbol z
	+\boldsymbol\phi\cdot \nabla \boldsymbol z
	-\left( \mathrm{div}\, \boldsymbol z\right)\boldsymbol\phi-\boldsymbol z\cdot \nabla \boldsymbol\phi\right)\nonumber\\
&= \left(\sigma(\boldsymbol y),\boldsymbol\phi\cdot \nabla \boldsymbol z
	-\boldsymbol z\cdot \nabla \boldsymbol\phi\right)\nonumber\\
	&=b\left(\phi,z, \sigma(y)\right) 
	-b\left(z,\phi,\sigma(y)\right).\end{align}
In view of this result, the state equation is to be understood in the sense of the following definition.
\begin{definition} Let $u\in L^2(I;H({\rm curl};\Omega))$ and
 $y_0\in W\cap H^3(\Omega)$. A function $ y\in L^\infty(I;  W\cap H^3(\Omega))$ with
$ \frac{\partial y }{\partial t}\in L^\infty(I; V)$
 is a solution of $(\ref{equation_etat})$ if $ y(0)= y_0$ in $\Omega$ and
	\begin{align}\label{var_form_state}
	&\left(\tfrac{\partial y(t) }{\partial t},\phi\right)
	+	2\alpha \left(D\tfrac{\partial y(t)}{\partial t}, D\phi\right)
	+2\nu\left(Dy(t),D\phi\right)
		+b\left( \phi, y(t), \sigma\left(y(t)\right)\right)-b\left(y(t),
	 \phi, \sigma\left(y(t)\right)\right)\nonumber\\
			&=\left( u(t),\phi\right)\qquad \mbox{for all} \ 
	\phi\in  V.\end{align}
\end{definition}
The result in the next proposition deals with the solvability of (\ref{equation_etat}) and is proved in \cite{BR03}. For the convenience of the reader, the corresponding estimates are derived hereafter. 
\begin{proposition}
 \label{existence_state} Assume that $u\in L^2(I;H({\rm curl};\Omega))$ and that $ y_0\in W\cap H^3(\Omega)$. Then the problem $(\ref{equation_etat})$ admits a unique solution 
  $ y\in L^\infty\left(I; W\cap H^3(\Omega)\right)$ with
$\displaystyle \tfrac{\partial y}{\partial t}\in L^2(I;V).$ Furthermore, the following estimates hold
\begin{align}
\label{state_est012} 
\left\|  y\right\|_{L^\infty(I;L^2(\Omega))}^2+2\alpha
 \left\|Dy\right\|_{L^\infty(I;L^2)}^2
+2\left\|Dy\right\|_{2,Q}^2\leq 4\left(
\left\|  y_0\right\|_2^2+2\alpha
 \left\|Dy_0\right\|_2^2+ \|u\|^2_{L^1(I;L^2)}\right),
\end{align}
\begin{align}
\label{state_est3}
\left\|{\rm curl}\,\sigma(y)\right\|_{L^\infty(I;L^2)}^2
&\leq 
2\left\|{\rm curl}\,\sigma(y_0)\right\|_2^2+
2\|{\rm curl} \, y\|_{L^\infty(I;L^2)}^2+4\left\|
	{\rm curl} \, u\right\|_{L^1(0,T;L^2)}^2,
	\end{align}
\begin{align}
\label{state_est4}
\alpha \left\|y\right\|_{L^\infty{(I;H^3)}}&\leq 
c\left(\left\|y\right\|_{L^\infty{(I;H^1)}}+\left\|{\rm curl}\,\sigma(y)\right\|_{L^\infty(I;L^2)}\right),	
\end{align}
\begin{align}
\label{state_t}
\left\|\tfrac{\partial y }{\partial t}\right\|_{2,Q}^2+
	\alpha \left\|D\tfrac{\partial y }{\partial t}\right\|_{2,Q}^2\leq c\left(\left\|u\right\|_{2,Q}^2
	+\left(\nu^2
	+\left\| {\rm curl}\,\sigma(y)\right\|_{L^\infty(I;L^2)}^2\right)
	\left\|y\right\|_{L^2(I;H^2)}^2\right),
\end{align} 
where $c$ is a positive constant only depending on $\Omega$.
\end{proposition}
{\bf Proof.} The proof is split into four steps.
\vspace{2mm}\\
{\it Step 1. $H^1$ in space estimate for $y$}.
By setting $\phi=y(t)$ in the variational formulation (\ref{var_form_state}), we obtain
 \begin{align}
	\tfrac{1}{2} \tfrac{d }{d t}\left(
	\left\|  y(t)\right\|_2^2+2\alpha
	 \left\|Dy(t)\right\|_2^2\right)
	 &+2\nu\left\|Dy(t)\right\|_2^2
	 \leq \left\|u(t)\right\|_2\left\|y(t)\right\|_2\nonumber\\
	 	&\leq \left\|u(t)\right\|_2
	 	\left(\left\|  y(t)\right\|_2^2+2\alpha
	 \left\|Dy(t)\right\|_2^2\right)^{\frac{1}{2}}
	 \label{en_y}
	\end{align}
yielding
$$
\tfrac{d }{d t}\left(\left\|  y(t)\right\|_2^2+2\alpha
	 \left\|Dy(t)\right\|_2^2\right)^{\frac{1}{2}}
\leq  \left\|u(t)\right\|_2.$$
Upon integration, we obtain
	$$\left(\left\|y(t)\right\|_2^2+2\alpha \left\|Dy(t)\right\|_{2}^2\right)^\frac{1}{2}
	\leq \left(\left\|y_0\right\|_2^2+2\alpha \left\|Dy_0\right\|_{2}^2\right)^\frac{1}{2}
	+\|u\|_{L^1\left(I;L^2\right)}$$
which, together with (\ref{en_y}), implies that
	$$\begin{array}{ll}\left\|y(t)\right\|_2^2+2\alpha \left\|Dy(t)\right\|_{2}^2+2\nu\left\|Dy(t)\right\|_{2}^2& \displaystyle 
	\leq \left\|y_0\right\|_2^2+2\alpha \left\|Dy_0\right\|_{2}^2+\int_0^t \left\|u(s)\right\|_2\left\|y(s)\right\|_{2}\,ds\vspace{2mm}\\
	&\leq \left\|y_0\right\|_2^2+2\alpha \left\|Dy_0\right\|_{2}^2+\|u\|_{L^1\left(I;L^2\right)} \left\|y\right\|_{L^\infty\left(I;L^2\right)}\vspace{2mm}\\
	&\leq 4\left(\left\|y_0\right\|_2^2+2\alpha \left\|Dy_0\right\|_{2}^2+\|u\|_{L^1\left(I;L^2\right)}^2\right).
	\end{array}$$
{\it Step 2. $L^2$ estimate in space for  ${\rm curl}\,\sigma(y)$.} 
	 Applying the curl to 
$(\ref{equation_etat})$, we obtain
	$$\tfrac{\partial}{\partial t}\left({\rm curl}\,\sigma(y(t))\right)-\nu \Delta \left({\rm curl} \, y(t)\right)+y(t)\cdot \nabla\left({\rm curl} \,\sigma(y(t))\right)={\rm curl} \, u(t)$$
which, for $\alpha\neq 0$, is equivalent to
	$$\tfrac{\partial}{\partial t}\left({\rm curl}\,\sigma(y(t))\right)
	+\tfrac{\nu}{\alpha}\, {\rm curl}\,\sigma(y(t))
	+y(t)\cdot \nabla\left({\rm curl} \,\sigma(y(t))\right)=
	{\rm curl} \, u(t)+\tfrac{\nu}{\alpha}\,{\rm curl} \, y(t)$$
Taking the $L^2$ scalar product with ${\rm curl} \, \sigma(y(t))$ yields
	$$\begin{array}{ll}
	&\tfrac{1}{2}\,\tfrac{d}{dt}\left\|{\rm curl}\,\sigma(y(t))
	\right\|_2^2+
	\tfrac{\nu}{\alpha}\left\|\mathrm{curl}\, 
	\sigma(y(t))\right\|_{2}^2=
	\left({\rm curl} \, u(t)+\tfrac{\nu}{\alpha}\,{\rm curl} \, y(t),
{\rm curl} \, \sigma(y(t))\right)\vspace{2mm}\\
	&\leq\left\|\mathrm{curl} \, u(t)\right\|_2\left\|{\rm curl} \, \sigma(y(t))\right\|_{2}
	+\tfrac{\nu}{2\alpha}\left\|\mathrm{curl} \, y(t)\right\|_2^2
	+\tfrac{\nu}{2\alpha}\left\|{\rm curl} \, \sigma(y(t))\right\|_{2}^2\end{array}$$
and thus
	$$\begin{array}{ll}\tfrac{d}{dt}\left\|{\rm curl}\,\sigma(y(t))
	\right\|_2^2+
	\tfrac{\nu}{\alpha}\left\|\mathrm{curl}\, 
	\sigma(y(t))\right\|_{2}^2&=e^{-\frac{\nu t}{\alpha}}
	\tfrac{d}{dt}\left(e^{\frac{\nu t}{\alpha}}
	\left\|{\rm curl}\,\sigma(y(t))\right\|_2^2\right)\vspace{2mm}\\
	&\leq
	\tfrac{\nu}{\alpha}\left\|\mathrm{curl} \, y(t)\right\|_2^2+
	2\left\|\mathrm{curl} \, u(t)\right\|_2\left\|{\rm curl} \, \sigma(y(t))\right\|_{2}.\end{array}$$
Multiplying both sides by $\text{e}^{\frac{\nu t}{\alpha}}$ and integrating, we obtain we obtainwe obtainwe obtainwe obtainwe obtainwe obtain we obtain we obtain we obtain we obtain
	$$\begin{array}{ll}&\left\|{\rm curl}\,\sigma(y(t))
	\right\|_2^2\vspace{2mm}\\
	&\displaystyle\leq \left\|{\rm curl}\,\sigma(y_0)	
	\right\|_2^2+\tfrac{\nu}{\alpha}
	\int_0^t e^{\frac{\nu}{\alpha}(s-t)}\left\|\mathrm{curl} \, y(s)\right\|_2^2\,ds\displaystyle+\int_0^t 2e^{\frac{\nu}{\alpha}(s-t)}\left\|\mathrm{curl} \, u(s)\right\|_2\left\|{\rm curl} \, \sigma(y(s))\right\|_{2}\,ds\vspace{2mm}\\
&\displaystyle\leq \left\|{\rm curl}\,\sigma(y_0)	
	\right\|_2^2+\tfrac{\nu}{\alpha}
	\left\|\mathrm{curl} \, y\right\|_{L^\infty(I;L^2)}^2
	\int_0^t e^{\frac{\nu}{\alpha}(s-t)}\,ds\vspace{2mm}\\
	&\displaystyle+2\left\|{\rm curl} \, \sigma(y)\right\|_{L^\infty(I;L^2)}\int_0^t e^{\frac{\nu}{\alpha}(s-t)}\left\|\mathrm{curl} \, u(t)\right\|_2\,ds\vspace{2mm}\\
&\displaystyle\leq \left\|{\rm curl}\,\sigma(y_0)	
	\right\|_2^2+
	\left\|\mathrm{curl} \, y\right\|_{L^\infty(I;L^2)}^2+2\left\|{\rm curl} \, \sigma(y)\right\|_{L^\infty(I;L^2)}\left\|\mathrm{curl} \, u\right\|_{L^1(0,T;L^2)}\vspace{2mm}\\
&\displaystyle\leq \left\|{\rm curl}\,\sigma(y_0)	
	\right\|_2^2+
	\left\|\mathrm{curl} \, y\right\|_{L^\infty(I;L^2)}^2+\tfrac{1}{2}
	\left\|{\rm curl} \, \sigma(y)\right\|_{L^\infty(I;L^2)}^2
	+2\left\|\mathrm{curl} \, u\right\|_{L^1(0,T;L^2)}^2.\end{array}$$
Estimate (\ref{state_est3}) follows immedialtely. \vspace{2mm}\\
{\it Step 3. $H^3$ estimate in space for $y$. } 
Since $\mathrm{curl}\,\Delta y(t)\in L^2(\Omega)$
and $\mathrm{div} \left(\mathrm{curl}\,\Delta y(t)\right)=0$, there exists a unique vector-potential $\psi\in H^1(\Omega)$ such that
	$$\left\{\begin{array}{ll}\rm{curl}\,
	 \psi=\mathrm{curl}\,\Delta y(t) &\quad \mbox{in} \ \Omega, \vspace{2mm}\\
	\nabla \cdot \psi=0 & \quad\mbox{in} \ \Omega,\vspace{2mm}\\
	\psi\cdot n=0 & \quad\mbox{on} \ \Gamma\end{array}\right.$$
and 
	\begin{equation}\label{sigma_phi}
	\left\|\psi\right\|_{H^1}\leq 
	c\left\|\mathrm{curl}\,\Delta y(t)\right\|_2.\end{equation}
It follows that
	$${\rm curl} \left(\Delta y(t)-\psi\right)=0$$
and there exists $\pi\in L^2(\Omega)$ such that
	$$\Delta y(t)-\psi+\nabla \pi=0.$$
Hence $y$ is the solution of the Stokes system
	$$\Delta y(t)+\nabla\pi=\psi$$
and satisfies
	\begin{equation}\label{y_phi}\left\|y(t)\right\|_{H^3}
	\leq c\left\|\psi\right\|_{H^1}.\end{equation}
Combining (\ref{sigma_phi}) and (\ref{y_phi}), we obtain
$$ \left\|y(t)\right\|_{H^3}\leq 
	c\left\|\mathrm{curl}\,\Delta y(t)\right\|_{2},$$
and thus
$$
\|y(t)\|_{H^3}\leq \tfrac{c}{\alpha}\left(
\left\|{\rm curl}\,\sigma(y(t))\right\|_2+
\left\|{\rm curl}\,y(t)\right\|_2
\right).$$
The estimate follows directly.
 \vspace{2mm}\\
{\it Step 4.  Estimates for the time derivative.}  
By setting  $\phi=\tfrac{\partial y }{\partial t}(t)$ in the variational formulation (\ref{var_form_state}), we obtain  
	$$\begin{array}{ll}
	\left\|\tfrac{\partial y }{\partial t}(t)\right\|_2^2+
	2\alpha \left\|D\tfrac{\partial y }{\partial t}(t)\right\|_2^2&=
	\left(u(t)+\nu\Delta y(t),\tfrac{\partial y }{\partial t}(t)\right) -\left( \boldsymbol{ curl}\,\sigma(\boldsymbol y(t))
	\times \boldsymbol y(t),
	\tfrac{\partial \boldsymbol y }{\partial t}(t)\right)\vspace{2mm}\\
	&\leq \left(\left\|u(t)\right\|_2+\nu\left\|\Delta y(t)\right\|_2
	+\left\|{\rm curl}\,\sigma(y(t))\right\|_2
	\left\|y(t)\right\|_\infty\right)
	\left\|\tfrac{\partial y }{\partial t}(t)\right\|_2\end{array}$$
and thus
	$$\begin{array}{ll}
	\left\|\tfrac{\partial y }{\partial t}(t)\right\|_2^2+
	\alpha \left\|D\tfrac{\partial y }{\partial t}(t)\right\|_2^2&\leq 
	c\left(\left\|u(t)\right\|_2^2
	+\nu^2\left\|\Delta y(t)\right\|_2^2
	+\left\|{\rm curl}\,\sigma(y(t))\right\|_2^2
	\left\|y(t)\right\|_\infty^2\right)\vspace{2mm}\\
	&\leq c\left(\left\|u(t)\right\|_2^2
	+\nu^2\left\|y(t)\right\|_{H^2}^2
	+\left\| {\rm curl}\,\sigma(y(t))\right\|_2^2
	\left\|y(t)\right\|_{H^2}^2\right).
	\end{array}$$
The claimed result follows then by integrating the previous inequality.
$\hfill\Box$
\section{Linearized state equation and analysis of the control-to-state mapping}
\setcounter{equation}{0}
As well known, the linearized equation associated with the state equation  plays a key role in the derivation of the necessary optimality conditions. Its solution coincides with the derivative of the control-to-state mapping and is related to the adjoint state through a suitable Green formula.\vspace{1mm}\\
The aim of this section is to establish the existence of a unique solution for an auxiliary linear system and to derive useful corresponding  estimates . The solvability of the linearized equation, seen as a particular case, the Lipchitz continuity and the G\^ateaux differentiability of the control-to-state mapping are then deduced from these results.\vspace{1mm}\\
Let $y_1, y_2$ be in $L^\infty\left(I; W\cap H^3(\Omega)\right)$ and consider the following linear equation
\begin{equation}\label{linearized}\left\{
  \begin{array}{ll}
   \displaystyle\tfrac{\partial \sigma(\boldsymbol z)}{\partial t}
  	-\nu \Delta \boldsymbol z+\boldsymbol{curl}\,
	\sigma(\boldsymbol z)
	\times  \boldsymbol y_1+
	\boldsymbol{curl}\,\sigma(\boldsymbol y_2)\times \boldsymbol z+\nabla \pi= \boldsymbol w&\quad\mbox{in} \ Q,\vspace{2mm}\\
	\nabla\cdot \boldsymbol z=0&\quad\mbox{in} \ Q,\vspace{2mm}\\
      \boldsymbol z\cdot\boldsymbol n=0, \quad (\boldsymbol n\cdot D \boldsymbol z)\cdot
      \boldsymbol \tau=0 &\quad\mbox{on}\ \Sigma,
      \vspace{2mm}\\      
      \boldsymbol z(0)=0 &\quad\mbox{in}\ \Omega, 
  \end{array}
\right.
\end{equation}
where  $w\in L^2(Q)$. In analogy to (\ref{equation_etat}), and taking into account  (\ref{trilinear_state}), we first propose the following definition for a solution of 
$(\ref{linearized})$.
\begin{definition} A function $ z\in L^\infty(I;W)$ with  
 $\frac{\partial z}{\partial t}\in L^2(I;V)$ is a solution of $(\ref{linearized})$ if $z(0)=0$ and
\begin{align}
\label{form_var_lin_2}
\left(\tfrac{\partial z(t)}{\partial t}, \phi\right)&+2\alpha
\left(D\tfrac{\partial z(t)}{\partial t},
D\phi\right)+
2\nu\left(Dz(t),D\phi\right)+
b\left( \phi, y_1(t), \sigma(z(t)) \right)-b\left( y_1(t), \phi, \sigma(z(t)) \right)\vspace{2mm}\nonumber\\
&+
b\left( \phi, z(t), \sigma(y_2(t)) \right)-b\left( z(t), \phi, \sigma(y_2(t)) \right)=
\left( w(t), \phi\right) \qquad \mbox{for all} \  \phi\in  V.
\end{align}
\end{definition}
The special Galerkin basis used to study the state equation (\ref{equation_etat}) does not seem appropriate to study the solvability of both the linearized equation considered in this section and the adjoint state equation that will be considered in Section \ref{sec_adjoint}. 
Indeed, the corresponding technique decomposes the problem into a mixed parabolic-hyperbolic system, looking for $z$ as the solution of a Stokes-like system and for $\mathrm{curl}\, \sigma(z)$ as the solution of the following transport equation
	$$\begin{array}{ll}\tfrac{\partial}{\partial t}\left(\mathrm{curl}\, \sigma(z(t))
	\right)+\tfrac{\nu}{\alpha}\,
	\mathrm{curl}\, \sigma(z(t))+ \,y_1(t)\cdot \nabla 
	\left(\mathrm{curl}\, \sigma(z(t))\right)
	+z(t)\cdot \nabla 
	\left(\mathrm{curl}\, \sigma(y_2(t))\right)\vspace{2mm}\\
	=
	\,\mathrm{curl}\, w(t)+\tfrac{\nu}{\alpha}\,\mathrm{curl}\, z(t).
	\end{array}$$
To (formally) derive the $L^2$ in space estimate for $\mathrm{curl}\, \sigma(z)$, let us multiply the transport equation by $\mathrm{curl}\, \sigma(z(t))$
	$$\begin{array}{ll}\tfrac{1}{2}\,\tfrac{d}{dt}
	\left(\left\|\mathrm{curl}\, \sigma(z(t))\right\|_2^2\right)+
	\tfrac{\nu}{\alpha}\left\|\mathrm{curl}\, \sigma(z(t))
	\right\|_2^2&=\left(\mathrm{curl}\, w(t)
	+\tfrac{\nu}{\alpha}\,\mathrm{curl}\, z(t),
	\mathrm{curl}\, \sigma(z(t))\right)\vspace{2mm}\\
	&-\left(z(t)\cdot \nabla 
	\left(\mathrm{curl}\, \sigma(y_2(t))\right),
	\mathrm{curl}\, \sigma(z(t))\right).\end{array}$$
The first term on the right-hand side can be easily handled by using the a priori $H^1$ estimates established in a first step. The second term, more delicate, can be managed if we guarantee that the coefficient $\mathrm{curl}\,\sigma(y_2)$ is $H^2$ in space. In other words, if the state variable $y_2$ belongs to $H^5$ in space. Following the regularity results stated in Theorem 5.3 and Remark 5.4 in \cite{CG97}, this property would be available if the pair $(u_2,y_0)$ belongs to $L^2(I;H^2(\Omega))\times H^5(\Omega)$. These difficulties are aggravated in the case of the adjoint equation because of the operators involved in its definition. \vspace{1mm}\\
According to these observations, and to the fact that $H^2$ a priori estimates for the linearized state and the adjoint state are sufficient to carry out our analysis, we will construct our solution by using the basis of eigenfunctions of the Stokes operator designed by Clopeau {\it et al.} to study the inviscid limit of the solutions of the Navier-Stokes equation with  Navier-slip boundary conditions \cite{CMR98}. In order to deal with the pressure term, the standard way to obtain $H^2$ a priori estimates in space would be to (formally) multiply equation (\ref{linearized}) by $\mathbb{A}z$ and to integrate. In our case, the main difficulty is then related with the term
	$$\left(\boldsymbol{curl}\, \sigma(\boldsymbol z(t))\times \boldsymbol y_1(t), \mathbb{A}\boldsymbol z(t)\right)$$ 
and is overcome by taking advantage of the nice properties induced by the Navier-slip boundary conditions and stated in Lemma \ref{curl_trace} and Lemma \ref{sigma_psigma}.\vspace{1mm}\\
We first state a lemma that will be used to derive $H^1$ a priori estimates for the linearized state equation and the adjoint equation. 
Unlike the Dirichlet boundary conditions for which the corresponding proofs are straightforward, the Navier-slip boundary conditions are more delicate to handle and proving that the boundary terms, induced by the performed integrations by parts, are vanishing is not an obvious issue.
\begin{lemma}\label{rm2}
Let $ y \in  W\cap H^3(\Omega)$ and $ z\in  W$. Then
	\begin{equation}\label{rm2_lin}
	\left|\left( \boldsymbol{curl}\,  \sigma(\boldsymbol z)	\times  \boldsymbol y, \boldsymbol z\right)\right|
	\leq 
	c \left\|y\right\|_{H^3}
	\left((1+\alpha)
	\|z\|_2^2+\alpha \left\|D z\right\|_2^2
	\right),\end{equation}
where  $c$ is a positive constant only depending on $\Omega$. 
\end{lemma}
{\bf Proof.} By taking into account the identity (\ref{trilinear_state}), we obtain
$$\left( \boldsymbol{curl}\,\sigma(\boldsymbol z)\times  \boldsymbol y,\boldsymbol z\right)$$
	$$=b\left(z, y, \sigma(z)\right) 
	-b\left(y,z, \sigma(z)\right) 
	=b(z,y,z)-b(y,z,z)-\alpha 
	\left(z\cdot \nabla y-y\cdot \nabla z,\Delta z\right)$$
	$$=b(z,y,z)+\alpha 
	\left(\boldsymbol z\cdot \nabla \boldsymbol y-\boldsymbol y\cdot \nabla \boldsymbol z,\boldsymbol{curl}\left(\boldsymbol{curl}\, \boldsymbol z\right)\right)$$
	$$=b(z,y,z)+\alpha 
	\left(\boldsymbol{curl}\left(\boldsymbol z\cdot \nabla \boldsymbol y-\boldsymbol y\cdot \nabla \boldsymbol z\right),\boldsymbol{curl}\, \boldsymbol z\right)
	+\alpha I$$
	$$=b(z,y,z)
	+\alpha b\left(\boldsymbol z, \boldsymbol{curl}\, \boldsymbol y,\boldsymbol{curl}\, \boldsymbol z\right)
	 -\alpha b\left(\boldsymbol y, \boldsymbol{curl}\, \boldsymbol z,\boldsymbol{curl}\, \boldsymbol z\right)+2\alpha\sum_{k=1}^3\left(\nabla \boldsymbol z_k\times \nabla \boldsymbol y_k,\boldsymbol{curl}\,\boldsymbol z\right)+\alpha I$$
	\begin{equation}\label{trilin_0}=b(z,y,z)
	+\alpha b\left(\boldsymbol z, \boldsymbol{curl}\, \boldsymbol y,\boldsymbol{curl}\, \boldsymbol z\right)+2\alpha\sum_{k=1}^2\left(\nabla \boldsymbol z_k\times \nabla \boldsymbol y_k,\boldsymbol{curl}\,\boldsymbol z\right)+\alpha I,
	\end{equation}
where 
	$$I= \displaystyle\int_\Gamma
	\left(y\cdot \nabla z-z\cdot \nabla y\right)\cdot \tau  \left(z\cdot g\right)\,dS.$$
Extending the exterior normal $n$ (defined a priori only on the boundary $\Gamma$) inside $\Omega$ by a vector field still denoted by $n$, using the Green formula and standard calculation, we can prove that for every $w\in H^2(\Omega)$ we have
	$$\displaystyle\int_{\Gamma}   \left(z\cdot \nabla y-y\cdot \nabla z\right)\cdot w\,dS$$
	$$=\displaystyle\left(z\cdot \nabla y-y\cdot \nabla z,
	w\,\mathrm{div} \, n\right)+
	\int_\Omega n\cdot \nabla\left( 
	\left(z\cdot \nabla y-y\cdot \nabla z\right)\cdot w\right) \,dx	$$
	$$=\displaystyle\left(z\cdot \nabla y-y\cdot \nabla z,
	w\,\mathrm{div} \, n\right)+b\left(n,w,z\cdot \nabla y-y\cdot \nabla z\right)+\left(\nabla
	\left(z\cdot \nabla y-y\cdot \nabla z\right)n,w\right)$$
	$$=\displaystyle\left(z\cdot \nabla y-y\cdot \nabla z,
	w\,\mathrm{div} \, n\right)+b\left(n,w,z\cdot \nabla y-y\cdot \nabla z\right)+\left(\nabla
	\left(z\cdot \nabla y-y\cdot \nabla z\right),w\otimes n\right)
	$$
	$$=\displaystyle\left(z\cdot \nabla y-y\cdot \nabla z,
	w\,\mathrm{div} \, n\right)+b\left(n,w,z\cdot \nabla y-y\cdot \nabla z\right)$$
	$$+b\left(z,\nabla y,w\otimes n\right)-b\left(y,\nabla z,w\otimes n\right)+\left(\nabla z\nabla y-\nabla y\nabla z,w\otimes n\right)$$
	$$=\displaystyle b\left(z,y,w \,\mathrm{div} \, n\right)-
	b\left(y,z,w \,\mathrm{div} \, n\right)+b\left(n,w,z\cdot \nabla y-y\cdot \nabla z\right)$$
	\begin{equation}\label{green_convective}
	 -b\left(z,w\otimes n,\nabla y\right)+
	b\left(y,w\otimes n,\nabla z\right)+\left(\nabla z\nabla y-\nabla y\nabla z,w\otimes n\right),
	\end{equation}
where 
	$$w\otimes n=\left(w_in_j\right)_{ij} \qquad \mbox{and} \qquad 
	b(\phi,\eta,\zeta)=\displaystyle \sum_{i,j,k}\int_\Omega \phi_k\tfrac{\partial \eta_{ij}}{\partial x_k}\zeta_{ij}\,dx$$
 Combining (\ref{trilin_0}) and (\ref{green_convective}), we deduce that
	\begin{align}\label{trilin_1}
	\left( \boldsymbol{curl}\,\sigma(\boldsymbol z)\times  \boldsymbol y,\boldsymbol z\right)
	&=b(z,y,z)+\alpha b\left(\boldsymbol z, \boldsymbol{curl}\, \boldsymbol y,\boldsymbol{curl}\, \boldsymbol z\right)+2\displaystyle\alpha\sum_{k=1}^2\left(\nabla \boldsymbol z_k\times \nabla \boldsymbol y_k,\boldsymbol{curl}\,\boldsymbol z\right)\nonumber\\
	&-\displaystyle \alpha b\left(z,y,\left(z\cdot g\right)\tau \,\mathrm{div} \, n\right)+\alpha\,
	b\left(y,z,\left(z\cdot g\right)\tau \,\mathrm{div} \, n\right)\nonumber\\
	&-\alpha b\left(n,\left(z\cdot g\right)\tau,z\cdot \nabla y-y\cdot \nabla z\right)
	- \alpha b\left(y,\left(z\cdot g\right)\tau\otimes n,\nabla z\right)\nonumber\\
	&+\alpha b\left(z,\left(z\cdot g\right)\tau\otimes n,\nabla y\right)
	+\alpha\left(\nabla z\nabla y-\nabla y\nabla z,\left(z\cdot g\right)\tau\otimes n\right).
	\end{align}
Standard arguments together with the Sobolev inequality give
	$$\left|b(z,y,z)\right|+\alpha\left|b\left(\boldsymbol z, \boldsymbol{curl}\, \boldsymbol y,\boldsymbol{curl}\, \boldsymbol z\right)\right|+\displaystyle 2\alpha\sum_{k=1}^2\left|\left(\nabla \boldsymbol z_k\times \nabla \boldsymbol y_k,\boldsymbol{curl}\,\boldsymbol z\right)\right|$$
	$$\leq \|z\|_2^2 \left\|\nabla y\right\|_\infty+
	\alpha\left(\|z\|_4 \left\|\nabla \boldsymbol{curl}\,\boldsymbol y\right\|_{4} \left\|\boldsymbol{curl}\, \boldsymbol z\right\|_2+
	2\sum_{k=1}^2\left\|\nabla \boldsymbol z_k\right\|_2
	\left\|\nabla \boldsymbol y_k\right\|_\infty\left\|\boldsymbol{curl}\,
	\boldsymbol z\right\|_2\right)$$
	\begin{equation}\label{trilin_2}\leq c\left\|y \right\|_{H^3}
	 \left(
	 \left\|\boldsymbol z\right\|_2^2+ \alpha
	\|\nabla z\|_2^2
	\right)	\end{equation}	
and
	$$\left|b\left(z,y, \left(z\cdot g\right)\tau \,\mathrm{div} \, n\right)-b\left(y,z, \left(z\cdot g\right)\tau \,\mathrm{div} \, n\right)\right|$$
	$$\leq \left(\left\|z\right\|_4\left\|\nabla y\right\|_2+
	\left\|y\right\|_4\left\|\nabla z\right\|_2
	\right)\left\| \left(z\cdot g\right)\tau \,\mathrm{div} \, n\right\|_4 $$
	$$\leq \left(\left\|z\right\|_4\left\|\nabla y\right\|_2+
	\left\|y\right\|_4\left\|\nabla z\right\|_2
	\right)\left\|z\right\|_4\left\|g\right\|_\infty
	\left\|\tau\right\|_\infty \left\|\mathrm{div} \, n\right\|_\infty	$$
	\begin{equation}\label{trilin_3}\leq c\left\|n\right\|_\infty^2\left\|\nabla n\right\|_\infty^2\left\|\nabla y\right\|_2\left\|\nabla z\right\|_2^2\leq cc_1(n)^4\left\|y\right\|_{H^3}\left\|\nabla z\right\|_2^2
	 \end{equation}
with $c_k(n)=\|n\|_{C^k(\bar\Omega)}$ and $c$ only depending on $\Omega$. Similarly, we have
	$$\left| b\left(n,\left(z\cdot g\right)\tau,z\cdot \nabla y-y\cdot \nabla z\right)+b\left(y,\left(z\cdot g\right)\tau\otimes n,\nabla z\right)-b\left(z,\left(z\cdot g\right)\tau\otimes n,\nabla y\right)
\right|$$
$$\leq \left\|n\right\|_\infty
	\left(\left\|z\right\|_4\left\|\nabla y\right\|_4+
	\left\|y\right\|_\infty\left\|\nabla z\right\|_2\right)
	\left\|\nabla\left(\left(z\cdot g\right)\tau\right)\right\|_2$$
$$
	+\left(\left\|y\right\|_\infty\left\|\nabla z\right\|_2+
	\left\|z\right\|_4\left\|\nabla y\right\|_4 \right)\left\|\nabla\left( \left(z\cdot g\right)\tau\otimes n\right)\right\|_2
$$
$$\leq c\left\|y\right\|_{H^3}\left\|\nabla z\right\|_2 \left(\left\|n\right\|_\infty\left\|\nabla\left(\left(z\cdot g\right)\tau\right)\right\|_2+\left\|\nabla\left( \left(z\cdot g\right)\tau\otimes n\right)\right\|_2\right)
$$
$$\leq c\left\|y\right\|_{H^3}\left\|\nabla z\right\|_2 
	\left\|\nabla\left(z\cdot g\right)\right\|_2
	\left(\|n\|_\infty\|\tau\|_\infty+\left\|\tau\otimes n\right\|_\infty\right)$$
$$+c\left\|y\right\|_{H^3}\left\|\nabla z\right\|_2 
\left\|z\cdot g\right\|_2\left(\|n\|_\infty
	\left\|\nabla\tau\right\|_\infty+
	\left\|\nabla\left(\tau\otimes n\right)\right\|_\infty\right)
$$
$$\leq c\left\|y\right\|_{H^3}
	\left\|n\right\|_\infty^2\left( \left\|n\right\|_\infty\left\|\nabla n\right\|_\infty+\left\|\nabla n\right\|_\infty^2
	+\left\|n\right\|_\infty\left\|\nabla^{(2)} n\right\|_\infty \right)\left\|\nabla z\right\|_2^2$$
	\begin{equation}\label{trilin_4}\leq c c_2(n)^4 \left\|y\right\|_{H^3}\left\|\nabla z\right\|_2^2
		\end{equation}
and 
	$$\left|\left(\nabla z\nabla y-\nabla y\nabla z,\left(z\cdot g\right)\tau\otimes n\right)\right|$$
	$$\leq 2 \left\|\nabla z\right\|_2  \left\|\nabla y\right\|_\infty  \left\|\left(z\cdot g\right)\tau\otimes n\right\|_2\leq c\left\|n\right\|_\infty^3
	 \left\|\nabla n\right\|_\infty 
	\left\|\nabla y\right\|_\infty \left\|\nabla z\right\|_2
	\left\|z\right\|_2$$
	\begin{equation}\label{trilin5}\leq cc_1(n)^4
	\left\|y\right\|_{H^3} \left\|\nabla z\right\|_2^2.\end{equation}
The claimed result follows then by combining (\ref{trilin_1})-(\ref{trilin_4}) and using the Korn inequality.$\hfill\Box$\vspace{2mm}\\
Now we are able to deal with the solvability of problem (\ref{linearized}). This is the aim of the next result.
\begin{proposition} \label{ex_uniq_lin}
Let $w$ be in $L^2(Q)$ and let $y_1$, $y_2$ be in $L^\infty(I;W\cap H^3(\Omega))$. Then equation $(\ref{linearized})$ admits a unique solution $z\in L^\infty(I;W)$ with  $\frac{\partial z}{\partial t}\in L^2(I;V)$. Moreover, the following estimates hold
$$\left\|  z\right\|_{L^\infty(I;L^2)}^2+2\alpha\left\|Dz\right\|_{L^\infty(I;L^2)}^2\leq
\left\| w\right\|^2_{2,Q} \ \displaystyle
e^{cT \left(1+\left(1+\alpha \right)\left\|y_1\right\|_{L^\infty(I;H^3)}\right)},$$		
$$\left\|Dz\right\|_{L^\infty(I;L^2)}^2+\alpha 
\left\|\mathbb{A}z\right\|_{L^\infty(I;L^2)}^2\leq \left(\tfrac{1}{\nu}\left\| w\right\|^2_{2,Q}+cT(1+\alpha)M
\left\|z\right\|_{L^\infty(I;L^2)}^2\right) e^{\frac{c(1+\alpha)T}{\alpha}M},$$	
$$\left\|\tfrac{\partial z}{\partial t}\right\|_{2,Q}^2+
\alpha\left\|D\tfrac{\partial z}{\partial t}\right\|_{2,Q}^2\leq 
c\left(\left\|w\right\|_{2,Q}^2+T\left(\nu^2+\tfrac{(1+\alpha)^2}{\alpha}M^2\right)
	\left\|z\right\|_{L^\infty(I;H^2)}^2\right),$$
where $M=\left\|y_1\right\|_{L^\infty(I;H^3)}+\left\|y_2\right\|_{L^\infty(I;H^3)}$ and where $c$ is a constant only depending on $\Omega$.
\end{proposition}
{\bf Proof.} Following \cite{CMR98}, there exists a set of eigenfunctions $\left(e_j\right)_j\subset  H^3(\Omega)$ of the problem
	\begin{equation}\label{eigen_functions}
	\left\{\begin{array}{ll}-\Delta e_j+\nabla \pi_j=
	\lambda_je_j &\quad \mbox{in} \ \Omega,\vspace{2mm}\\
	\nabla\cdot e_j=0& \quad\mbox{in} \ \Omega,\vspace{2mm}\\
	e_j\cdot n=0, \qquad 
	\left(De_j\, n\right)\cdot \tau=0&\quad \mbox{on} \ \Gamma,\end{array}\right.
	\end{equation}
with 
	$$0<\lambda_1<\cdots<\lambda_j<\cdots \longrightarrow +\infty.$$
The functions $e_j$ form an orthonormal basis in $H$. The approximate problem is defined by
	\begin{equation}\label{faedo_galerkin_lin}
	\left\{\begin{array}{llll}\displaystyle \mbox{Find} \ z_k(t)
	=\sum_{i=1}^k \zeta_{i}(t) e_i \ \mbox{solution}, \ \mbox{for} \ 
	1\leq j\leq k, \ \mbox{of} \vspace{2mm}\\ 
	\displaystyle	\left(\tfrac{\partial z_k(t) }{\partial t}, e_j\right)+2\alpha
	\displaystyle\left(D\tfrac{\partial z_k(t) }{\partial t},
	De_j\right)+2\nu\left(Dz_k(t),De_j\right)\vspace{2mm}\\
	\qquad \qquad \,\quad  +\left(\boldsymbol{ curl}\,\sigma(\boldsymbol z_k(t))\times 
	\boldsymbol y_1(t)+
	 \boldsymbol{ curl}\,\sigma(\boldsymbol y_2(t))\times
	  \boldsymbol z_k(t), \boldsymbol e_j\right)=
	\left(w(t), e_j\right),\vspace{2mm}\\
	z_k(0)=0.
	\end{array}\right.
	\end{equation}
This is a linear differential system for $\zeta (t)=(\zeta_1 (t),\dots, \zeta_k (t) )^\top$, of the form
	$$
	A\tfrac{d\zeta(t)}{dt} +N(t)\zeta (t)=F(t),
	$$
where $A$ is the nonsingular constant matrix defined by
	$$A_{ij}=\left(e_i,e_j\right)+2\alpha\left(De_i,De_j\right) \qquad i,j=1,\cdots,k$$
 $N(t)$ is the matrix given by
	  $$\begin{array}{ll}N_{ij}(t)&=2\nu\left(De_i,De_j\right)+b\left(e_i,y_1(t),\sigma(e_j)\right)-b\left(y_1(t),e_i,\sigma(e_j)\right)\vspace{2mm}\\
	&+b\left(e_i,e_j,\sigma(y_2(t))\right)-
	b\left(e_j,e_i,\sigma(y_2(t))\right) \qquad i,j=1,\cdots,k\end{array}$$ 
and 
	$$F_i(t)=\left(w(t),e_i\right)\qquad i=1,\cdots,k.$$
Taking into account the regularity properties of $y$ and the assumptions on $w$, we deduce that $N\in L^\infty(0,T)$ and $F\in L^2(0,T)$ and, thus, the previous differential system admits a unique solution $\zeta\in H^1(0,T)$.\vspace{2mm}\\
The remaing proof is split into four steps.\vspace{2mm}\\
{\it Step 1. Uniform  $H^1$ in space estimate}.  Multiplying both sides of the equation (\ref{faedo_galerkin_lin}) by $\zeta_j (t)$ and taking the sum over $j=1,\dots,k$, we easily verify that 
  $ z_k$ satisfies the  equality:
\begin{align}
\label{P_zeta}
\displaystyle\tfrac{1}{2}\,\tfrac{d }{d t}\left(\left\|  z_k(t)\right\|_2^2+2\alpha
\left\|Dz_k(t)\right\|_2^2\right)&+2\nu \left\|D  z_k(t)\right\|_2^2
=\left( w(t),z_k(t)\right)-\left( \boldsymbol{ curl}\,\sigma(\boldsymbol z_k(t))\times \boldsymbol y_1(t),
\boldsymbol z_k(t)\right).\nonumber
\end{align}
By using (\ref{rm2_lin}), we  estimate the right hand side 
\begin{align}
&\tfrac{d }{d t}\left(\left\|  z_k(t)\right\|_2^2+2\alpha
\left\|Dz_k(t)\right\|_2^2\right)+4\nu \left\|Dz_k(t)\right\|_2^2\nonumber\\
&\leq  \left\| w(t)\right\|^2_2+\|z_k(t)\|_2^2+c\left\|y_1(t)\right\|_{H^3}\left((1+\alpha)
\|z_k(t)\|_2^2+\alpha \left\|D z_k(t)\right\|_2^2\right)\nonumber\\
&	\leq  \left\| w(t)\right\|^2_2+c \left(1+\left(1+\alpha \right)\left\|y_1(t)\right\|_{H^3}\right)
\left(\left\|z_k(t)\right\|_2^2+2\alpha \left\|Dz_k(t)\right\|_2^2\right)\nonumber\\
&	\leq  F_1(t)+G_1
\left(\left\|z_k(t)\right\|_2^2+2\alpha \left\|Dz_k(t)\right\|_2^2\right).
\end{align}
where 
	$$F_1(t)=\left\| w(t)\right\|^2_2 \qquad \mbox{and} \qquad
	G_1=c \left(1+(1+\alpha)\|y_1\|_{L^\infty(I;H^3)}\right).$$
Upon integration, we obtain
$$\begin{array}{ll}
&\displaystyle \left\|  z_k(t)\right\|_2^2+2\alpha
\left\|Dz_k(t)\right\|_2^2+\int_0^t 4\nu \left\|Dz_k(s)\right\|_2^2ds\vspace{1mm}\\
&\displaystyle\leq \left\|F_1\right\|_{L^1(0,t)}
+ G_1\int_0^t\left(\left\|  z_k(s)\right\|_2^2+2\alpha
\left\|Dz_k(s)\right\|_2^2\right)ds.
\end{array}$$	
Therefore, by using Gronwall's inequality we deduce that
\begin{equation}\label{est_zk_H1}
\left\|  z_k\right\|_{L^\infty(I;L^2)}^2+2\alpha
\left\|Dz_k\right\|_{L^\infty(I;L^2)}^2\leq
\left\|F_1\right\|_{L^1(0,T)}
\displaystyle
e^{G_1T}.
\end{equation}		
{\it Step 2. Uniform $H^2$ in space estimate}. 
Taking into account (\ref{eigen_functions}), we have $\mathbb{A} e_j=\lambda_j e_j$, and thus
	$$	\left(\tfrac{\partial z_k(t)}{\partial t},\lambda_je_j\right)=\displaystyle
	\sum_{i=1}^k \zeta_{i}'(t)\left(e_i,\lambda_je_j\right)=
	\sum_{i=1}^k \zeta_{i}'(t)\left(2De_i,De_j\right)=\left(2\tfrac{\partial }{\partial t}Dz_k(t),De_j\right).$$
Similarly, we have
	$$\begin{array}{ll}\left(2D \tfrac{\partial z_k(t)}{\partial t},\lambda_jD e_j\right)&=
	-\left(\Delta \tfrac{\partial z_k(t)}{\partial t},\lambda_je_j\right)=-
	\left(\Delta \tfrac{\partial z_k(t)}{\partial t},\mathbb{A} e_j\right)=\left(\tfrac{\partial}{\partial t}\Delta z_k(t),-\mathbb{A} e_j\right)\vspace{2mm}\\
	&=
	\left(\tfrac{\partial}{\partial t}\mathbb{A} z_k(t),\mathbb{A} e_j\right)\end{array}$$
and 
$$\left(2D z_k(t),\lambda_jD e_j\right)=
	-\left(\Delta z_k(t),\lambda_je_j\right)=
	-\left(\Delta z_k(t),\mathbb{A} e_j\right)=
	\left(\mathbb{A} z_k(t),\mathbb{A} e_j\right).$$
Multiplying (\ref{faedo_galerkin_lin}) by $\lambda_{j}$, we obtain
	$$\begin{array}{ll}\left(2\tfrac{\partial }{\partial t}Dz_k(t),De_j\right)+\alpha
	\left(\tfrac{\partial}{\partial t}\mathbb{A} z_k,\mathbb{A} e_j\right)+\nu
	\left(\mathbb{A} z_k(t),\mathbb{A} e_j\right)\vspace{3mm}\\
	=\left(w(t),\mathbb{A} e_j\right)-\left( \boldsymbol{curl}\, \sigma(\boldsymbol z_k(t))\times \boldsymbol y_1(t)+
	\boldsymbol{curl}\, \sigma(\boldsymbol y_2(t))\times
	  \boldsymbol z_k(t), \mathbb{A}e_j\right)\end{array}$$ 
yielding
\begin{align}
&\tfrac{1}{2}\tfrac{d}{d t}\left(2\left\|Dz_k(t)\right\|_2^2+\alpha \left\|\mathbb{A}z_k(t)\right\|_2^2\right)
+\nu \left\|\mathbb{A}z_k(t)\right\|_2^2\nonumber\\
&=\left( w(t), \mathbb{A}z_k(t)\right)-	\left( \boldsymbol{curl}\, \sigma(\boldsymbol z_k(t))\times  \boldsymbol y_1(t)+
	\boldsymbol{curl}\, \sigma(\boldsymbol y_2(t))\times
	\boldsymbol z_k(t), \mathbb{A}\boldsymbol z_k(t)\right).\nonumber
	\end{align}
Hence, by using the Young inequality 
\begin{align}
\tfrac{d}{d t}\left(2\left\|Dz_k(t)\right\|_2^2+\alpha \left\|\mathbb{A}z_k(t)\right\|_2^2\right)+\nu \left\|\mathbb{A}z_k(t)\right\|_2^2&\leq
\tfrac{1}{\nu}\left\| w(t)\right\|_2^2+2\left|\left( \boldsymbol{curl}\, \sigma(\boldsymbol z_k(t))\times \boldsymbol y_1(t), \mathbb{A}\boldsymbol z_k(t)\right)\right| \nonumber \\
&+2\left|\left(\boldsymbol{curl}\, \sigma(\boldsymbol y_2(t))\times
 \boldsymbol z_k(t), \mathbb{A}\boldsymbol z_k(t)\right)\right|.
\label{fesh2}
\end{align}
On one hand, we have
\begin{align}
&\left|\left(\boldsymbol{curl}\, \sigma(\boldsymbol z_k(t))\times \boldsymbol y_1(t), \mathbb{A}\boldsymbol z_k(t)\right)\right|\nonumber\\
&\leq\left|\left( \boldsymbol{curl}\, \sigma(\boldsymbol z_k(t))\times \boldsymbol y_1(t), \mathbb{A}\boldsymbol z_k(t)+\Delta \boldsymbol z_k(t)\right)\right|+\left|\left( \boldsymbol{curl}\, \sigma(\boldsymbol z_k(t))\times \boldsymbol y_1(t), 
	\Delta \boldsymbol z_k(t)\right)\right|\nonumber\\
&\leq\left|b\left( \mathbb{A}z_k(t)+\Delta z_k(t),y_1(t),\sigma(z_k(t))\right)-b\left(y_1(t), \mathbb{A}z_k(t)+\Delta z_k(t),\sigma(z_k(t))\right)\right|
\nonumber\\
&+\left|b\left(\Delta z_k(t),y_1(t),\sigma(z_k(t))\right)-b\left(y_1(t),\Delta z_k(t),\sigma(z_k(t))\right)\right|  \nonumber \\
&\leq\left|b\left( \mathbb{A}z_k(t)+\Delta z_k(t),y_1(t),\sigma(z_k(t))\right)-b\left(y_1(t), \mathbb{A}z_k(t)+\Delta z_k(t),\sigma(z_k(t))\right)\right|
\nonumber\\
&+\left|b\left(\Delta z_k(t),y_1(t),\sigma(z_k(t))\right)+b\left(y_1(t),z_k(t),\Delta z_k(t)\right)\right|  \nonumber \\
&\leq \left\|y_1(t)\right\|_{1,\infty}
\left(2\left\|\mathbb{A}z_k(t)+\Delta  z_k(t)\right\|_{H^1}
\left\|\sigma(z_k(t))\right\|_{2}+
\left\|\Delta z_k(t)\right\|_{2}\left(\left\|\sigma(z_k(t))\right\|_{2}+\left\| z_k(t)\right\|_{H^1}\right)\right) \nonumber 
	\end{align}
and applying (\ref{sigma_psigma2}) and (\ref{yh_sigma}), we deduce that
	\begin{align}
	\left|\left( \boldsymbol{curl}\, \sigma(\boldsymbol z_k(t))\times \boldsymbol y_1(t), \mathbb{A}\boldsymbol z_k(t)\right)\right|&\leq c(1+\alpha)\left\|y_1(t)\right\|_{1,\infty}
\left\|z_k(t)\right\|_{H^2}^2\nonumber\\
	&\leq  c(1+\alpha)\left\|y_1(t)\right\|_{H^3}
	\left(\left\|z_k(t)\right\|^2_{2}+\left\|\mathbb{A}z_k(t)\right\|^2_{2}
\right)\nonumber\\
	&\leq  c(1+\alpha)\left\|y_1\right\|_{L^\infty(I;H^3)}\left(\left\|z_k\right\|^2_{L^\infty(I;L^2)}+\left\|\mathbb{A}z_k(t)\right\|^2_{2}
\right).
\label{vis3}
	\end{align}
On the other hand, by using similar arguments we obtain	
\begin{align}
\left|\left(
\boldsymbol{curl}\, \sigma(\boldsymbol y_2(t))\times
\boldsymbol z_k(t), \mathbb{A}\boldsymbol z_k(t)\right)\right|&\leq 
\left\|{\rm curl}\, \sigma(y_2(t))\right\|_{2}
\left\|z_k(t)\right\|_{\infty}\left\|\mathbb{A}z_k(t)\right\|_{2}
\nonumber\\
&\leq 
c\left\|{\rm curl}\, \sigma(y_2(t))\right\|_{2}
\left\|z_k(t)\right\|^2_{H^2}\nonumber\\
&\leq 
c\left\|{\rm curl}\, \sigma(y_2(t))\right\|_{2}
\left(\left\|z_k(t)\right\|^2_{2}+\left\|\mathbb{A}z_k(t)\right\|^2_{2}
\right)\nonumber\\
	&\leq  c(1+\alpha)\left\|y_2(t)\right\|_{H^3}\left(\left\|z_k(t)\right\|^2_{2}+\left\|\mathbb{A}z_k(t)\right\|^2_{2}\right)\nonumber\\
	&\leq  c(1+\alpha)\left\|y_2\right\|_{L^\infty(I;H^3)}\left(\left\|z_k\right\|^2_{L^\infty(I;L^2)}+\left\|\mathbb{A}z_k(t)\right\|^2_{2}
\right).
\label{vis1}
	\end{align}
Combining (\ref{fesh2})-(\ref{vis1}), we obtain
\begin{align}
&\tfrac{d}{d t}\left(2\left\|Dz_k(t)\right\|_2^2+\alpha 
\left\|\mathbb{A}z_k(t)\right\|_2^2\right)+\nu
\left\|\mathbb{A}z_k(t)\right\|_2^2\nonumber\\
&\leq F_2(t)+G_2
\left(2\left\|Dz_k(t)\right\|_2^2+\alpha 
\left\|\mathbb{A}z_k(t)\right\|_2^2\right),
		\nonumber
		\end{align}
where
\begin{align}
&F_2(t)= \tfrac{1}{\nu}\left\| w(t)\right\|_2^2+c(1+\alpha)\left(\left\|y_1\right\|_{L^\infty(I;H^3)}+\left\|y_2\right\|_{L^\infty(I;H^3)}\right)
\left\|z_k\right\|^2_{L^\infty(I;L^2)},\nonumber\\
&G_2=\tfrac{c(1+\alpha)}{\alpha}\left(\left\|y_1\right\|_{L^\infty(I;H^3)}+\left\|y_2\right\|_{L^\infty(I;H^3)}\right).
		\nonumber
		\end{align}
Upon integration, we have
	$$2\left\|Dz_k(t)\right\|_2^2+\alpha 
\left\|\mathbb{A}z_k(t)\right\|_2^2\leq 
	\left\|F_2\right\|_{L^1(0,t)}+G_2\int_0^t 
\left(2\left\|Dz_k(s)\right\|_2^2+\alpha 
\left\|\mathbb{A}z_k(s)\right\|_2^2\right)\,ds$$
and by using Gronwall's inequality, we deduce that
\begin{align}
\label{lipschitz_est_V2_zm}
2\left\|Dz_k\right\|_{L^\infty(I;L^2)}^2+\alpha 
\left\|\mathbb{A}z_k\right\|_{L^\infty(I;L^2)}^2
&\leq \left\|F_2\right\|_{L^1(0,T)}\,e^{G_2T}.
\end{align}		
	{\it Step 3. Uniform estimate for the time derivative}. Let us multiply both sides of (\ref{faedo_galerkin_lin}) 
	by  $\frac{d\zeta_{j}(t)}{dt}$ and sum over $j=1, \dots,k$. This gives
	$$\begin{array}{ll}
	\left\|\tfrac{\partial z_k(t)}{\partial t}\right\|_2^2+
	2\alpha \left\|D\tfrac{\partial z_k(t)}{\partial t}\right\|_2^2&=
	\left(w(t)+\nu\, \Delta z_k(t),\tfrac{\partial z_k(t)}{\partial t}\right)\\
	&-\left( \boldsymbol{ curl}\,\sigma(\boldsymbol z_k(t))\times
	 \boldsymbol y_1(t)+
	\boldsymbol{ curl}\,\sigma(\boldsymbol y_2(t))\times
	  \boldsymbol z_k(t),\tfrac{\partial \boldsymbol z_k(t)}{\partial t}\right).
	\end{array}$$
Since $$\begin{array}{ll}
	&\left|	\left(\boldsymbol{ curl}\,\sigma(\boldsymbol z_k(t))\times \boldsymbol y_1(t)+
	\boldsymbol{curl}\,\sigma(\boldsymbol y_2(t))\times
	  \boldsymbol z_k(t),\tfrac{\partial \boldsymbol z_k(t)}{\partial t} \right)\right|\vspace{2mm}\\
	&\leq 
	  \left|b\left(\tfrac{\partial z_k(t)}{\partial t},y_1(t),\sigma(z_k(t))\right)   \right|
	+\left|b\left(y_1(t),\tfrac{\partial z_k(t)}{\partial t},
	\sigma(z_k(t))\right)   \right|\vspace{2mm}\\
	  &+\left|b\left(\tfrac{\partial z_k(t)}{\partial t},z_k(t),
	\sigma(y_2(t))\right)   \right|+\left|b
	  \left(z_k(t),\tfrac{\partial z_k(t)}{\partial t},
	\sigma(y_2(t))\right)   \right|\vspace{1mm}\\
	  &\leq c(1+\alpha)\left(\left\|y_1(t)\right\|_{H^3}+\left\|y_2(t)\right\|_{H^3}\right)\left\|z_k(t)\right\|_{H^2}\left\|\tfrac{\partial z_k(t)}{\partial t}\right\|_{H^1}\vspace{1mm}\\
	  &\leq c(1+\alpha)\left(\left\|y_1(t)\right\|_{H^3}+\left\|y_2(t)\right\|_{H^3}\right)\left\|z_k(t)\right\|_{H^2}\left(\left\|\tfrac{\partial z_k(t)}{\partial t}\right\|_{2}+\left\|D\tfrac{\partial z_k(t)}{\partial t}\right\|_{2}\right),
	\end{array}$$
 we deduce that
	$$\begin{array}{ll}
	\left\|\tfrac{\partial z_k(t) }{\partial t}\right\|_2^2+
	2\alpha \left\|D\tfrac{\partial z_k(t)}{\partial t}\right\|_2^2&\leq
	\left(\left\|w(t)\right\|_2+\nu \left\|\Delta z_k(t)\right\|_2\right)
	\left\|\tfrac{\partial z_k(t)}{\partial t}\right\|_2
	\vspace{1mm}\\
	&+c(1+\alpha)\left(\left\|y_1(t)\right\|_{H^3}+\left\|y_2(t)\right\|_{H^3}\right)
		\left\|z_k(t)\right\|_{H^2}\left\|\tfrac{\partial z_k(t)}{\partial t}\right\|_2\vspace{1mm}\\
	&+c(1+\alpha)\left(\left\|y_1(t)\right\|_{H^3}+\left\|y_2(t)\right\|_{H^3}\right)
		\left\|z_k(t)\right\|_{H^2}\left\|D\tfrac{\partial z_k(t)}{\partial t}\right\|_{2}
	\end{array}$$
and by using the Young inequality, we obtain
	$$\begin{array}{ll}
	\left\|\tfrac{\partial z_k(t) }{\partial t}\right\|_2^2+
	\alpha \left\|D\tfrac{\partial z_k(t)}{\partial t}\right\|_2^2\vspace{2mm}\\
	\leq
	c\left(\left\|w(t)\right\|_2^2+\left(\nu^2+\tfrac{(1+\alpha)^2}{\alpha}\left(\left\|y_1(t)\right\|_{H^3}+
	\left\|y_2(t)\right\|_{H^3}\right)^2\right)
	\left\|z_k(t)\right\|_{H^2}^2\right)\vspace{2mm}\\
	\leq
	c\left(\left\|w(t)\right\|_2^2+\left(\nu^2+\tfrac{(1+\alpha)^2}{\alpha}\left(\left\|y_1\right\|_{L^\infty(I;H^3)}+\left\|y_2\right\|_{L^\infty(I;H^3)}\right)^2\right)
	\left\|z_k\right\|_{L^\infty(I;H^2)}^2\right).
	\end{array}$$
Therefore 
\begin{align}
&\left\|\tfrac{\partial z_k}{\partial t}\right\|_{2,Q}^2+
\alpha\left\|D\tfrac{\partial z_k}{\partial t}\right\|_{2,Q}^2\nonumber\\
&\leq 
c\left(\left\|w\right\|_{2,Q}^2+T\left(\nu^2+\tfrac{(1+\alpha)^2}{\alpha}\left(\left\|y_1\right\|_{L^\infty(I;H^3)}+\left\|y_2\right\|_{L^\infty(I;H^3)}\right)^2\right)
	\left\|z_k\right\|_{L^\infty(I;H^2)}^2\right).
\label{lipschitz_est_V2_zktt}
\end{align}
{\it Step 4. Passage to the limit}.
According to (\ref{est_zk_H1}),  (\ref{lipschitz_est_V2_zm}) and  (\ref{lipschitz_est_V2_zktt}), the sequences $\left(z_k\right)_k$ and $\left(\tfrac{\partial z_k }{\partial t}\right)_k$ are uniformly bounded in $L^\infty(I;W)$ and $L^2(I;V)$, respectively. Then there exists a subsequence, still indexed by $k$, and function $ z\in  L^\infty(I;W)$ such that
	$$ z_k \longrightarrow  z \qquad \mbox{weakly* in} \ L^\infty(I;H^2(\Omega)),$$
	$$ \tfrac{\partial z_k }{\partial t}\longrightarrow  \tfrac{\partial z}{\partial t} \qquad \mbox{weakly in} \ L^2(I;H^1(\Omega)).$$
	These convergence results imply that $z\in C(\bar I; H^1(\Omega))$ and as 
	$z_k(0)$ converges to zero in $H^1(\Omega)$, we deduce that $z$ verifies the initial condition
	$z(0)=0.$
By passing to the limit in $(\ref{faedo_galerkin_lin})_1$, we obtain for every $j\geq 1$
$$\begin{array}{ll}
	\left(\tfrac{\partial z(t) }{\partial t}, e_j\right)&+2\alpha
	\left(D\tfrac{\partial z(t)}{\partial t},
	D e_j\right)+2\nu\left(Dz(t),D e_j\right)+
b\left(e_j, y_1(t), \sigma(z(t)) \right)-b\left( y_1(t),e_j, \sigma(z(t)) \right)\vspace{2mm}\\
	&+
b\left(e_j, z(t), \sigma(y_2(t)) \right)-b\left( z(t),e_j, \sigma(y_2(t)) \right) =
	\left(w(t),  e_j\right)
	\end{array}$$
and by density, we prove that $ z$ satisfies (\ref{form_var_lin_2}). Moreover, $z$ satisfies the claimed estimates.
Since (\ref{linearized}) is linear, the uniqueness result is a direct consequence of the first estimate.\hfill\Box\vspace{2mm}\\
Next, we derive some useful estimates related with the Lipschitz continuity of the state with respect to the control variable. More precisely, if $u_1$, $u_2$ are two controls
and if $y_{u_1}$, $y_{u_2}$ are the corresponding states then we are interested in estimating the difference
$y_{u_2}-y_{u_1}$ with respect to $u_2-u_1$ in adequate topologies. The result is a direct consequence of Proposition \ref{ex_uniq_lin}.
\begin{proposition}
\label{Lipschitz}
Let $ u_1,  u_2\in  L^2(I;H({\rm curl};\Omega))$ and let $y_{u_1},  y_{u_2}$
be the corresponding solutions of $(\ref{equation_etat})$. Then the following estimates hold \vspace{-2mm}
$$\left\| y_{u_2}-y_{u_1}\right\|_{L^\infty(I;L^2)}^2+2\alpha
\left\|D\left(y_{u_2}-y_{u_1}\right)\right\|_{L^\infty(I;L^2)}^2
\leq
\left\| u_2-u_1\right\|^2_{2,Q}
e^{cT \left(1+\left(1+\alpha \right)\left\|y_{u_1}\right\|_{L^\infty(I;H^3)}\right)},$$
$$\left\|D\left(y_{u_2}-y_{u_1}\right)\right\|_{L^\infty(I;L^2)}^2+\alpha 
\left\|\mathbb{A}\left(y_{u_2}-y_{u_1}\right)
	\right\|_{L^\infty(I;L^2)}^2$$
$$\leq \left(\tfrac{1}{\nu}\left\|  u_2-u_1\right\|^2_{2,Q}+cT(1+\alpha)M
\left\|y_{u_2}-y_{u_1}\right\|^2_{L^\infty(I;L^2)}\right)
	e^{\frac{c(1+\alpha)T}{\alpha}M}$$
where $M=\left\|y_{u_1}\right\|_{L^\infty(I;H^3)}+\left\|y_{u_2}\right\|_{L^\infty(I;H^3)}$ and where $c$ is a positive constant only depending on $\Omega$.
\end{proposition}
{\bf Proof.} It is easy to see that $y=y_{u_2}-y_{u_1}$ is the unique solution of the linear problem (\ref{linearized}) for $y_1=y_{u_1}$, $y_2=y_{u_2}$ and $w=u_2-u_1$. The claimed estimates are then a direct consequence of Proposition \ref{ex_uniq_lin}.$\hfill\Box$\vspace{2mm}\\
Finally, we are able to prove that the control-to-state mapping $u \mapsto y_u $ is G\^ateaux differentiable.
\begin{proposition}
\label{Gateau}
Let $\rho$ be such that $0<\rho<1$, $y_0\in  W\cap H^3(\Omega)$ and $(v,w)\in  \left(L^2(I;H({\rm curl};\Omega))\right)^2$. Set $ u_\rho= u+\rho  w$ and let $y_u$ and $ y_{\rho}$ be the solutions of 
$(\ref{equation_etat})$ corresponding to $u$ and  $ u_\rho$, respectively. Then, we have
	\begin{equation}
	\label{gateau_1}
	 y_{\rho}= y_u+\rho 
	 z_{uw}+\rho r_\rho \qquad 
	\mbox{with} \ \lim_{\rho\rightarrow 0}
	\left\|r_\rho\right\|_{L^\infty(I;H^1)}=0,
	\end{equation}
	and
	\begin{equation}
	\label{gateau_2}J\left(u_\rho,y_\rho\right)=J\left(u,y_u\right)+\rho 
	\int_0^T\left(\left(z_{uw},y_u-y_d\right)+\lambda(u,w)\right)dt
	+o(\rho),\end{equation}
	where $ z_{uw}\in  L^\infty(I;W)$ is the solution of the linearized equation
\begin{equation}\label{linearized_gateau}\left\{
  \begin{array}{ll}
   \displaystyle\tfrac{\partial \sigma(\boldsymbol z)}{\partial t}
  	-\nu \Delta \boldsymbol z+\boldsymbol{curl}\,
	\sigma(\boldsymbol z)
	\times  \boldsymbol {y_u}+
	\boldsymbol{curl}\,\sigma(\boldsymbol {y_u})\times \boldsymbol z+\nabla \pi= \boldsymbol w&\quad\mbox{in} \ Q,\vspace{2mm}\\
	\nabla\cdot \boldsymbol z=0&\quad\mbox{in} \ Q,\vspace{2mm}\\
      \boldsymbol z\cdot\boldsymbol n=0, \quad (\boldsymbol n\cdot D \boldsymbol z)\cdot
      \boldsymbol \tau=0 &\quad\mbox{on}\ \Sigma,
      \vspace{2mm}\\      
      \boldsymbol z(0)=0 &\quad\mbox{in}\ \Omega.
  \end{array}
\right.
\end{equation}
\end{proposition}
{\bf Proof.} Let us first notice that due to Proposition \ref{ex_uniq_lin}, equation (\ref{linearized_gateau}) admits a unique solution $z_{uw}$. Moreover, easy calculation shows that $z_{\rho}=\tfrac{ y_{\rho}- y_u}{\rho}$ satisfies
	$$\displaystyle\tfrac{\partial  \sigma\left(\boldsymbol z_\rho\right)}{\partial t}-\nu \Delta \boldsymbol  z_\rho+ \boldsymbol{ curl}\,
	 \sigma\left(\boldsymbol z_\rho\right)	\times \boldsymbol {y_u}+ \boldsymbol{curl}\,\sigma\left(\boldsymbol y_\rho\right)\times \boldsymbol z_\rho+\nabla\pi_\rho=\boldsymbol w,$$ 
and thus $r_\rho= z_\rho- z_{uw}$ satisfies
	$$
	\displaystyle\tfrac{\partial  \sigma\left(\boldsymbol r_\rho\right)}{\partial t}	
	-\nu \Delta \boldsymbol r_\rho+ 
	\boldsymbol{ curl}\,\sigma\left(\boldsymbol r_\rho\right)\times \boldsymbol {y_{u}}
	+\boldsymbol{ curl}\,\sigma\left(\boldsymbol y_\rho\right)\times\boldsymbol r_\rho
	+ \boldsymbol{ curl}\,\sigma\left(\boldsymbol y_\rho-\boldsymbol {y_u}\right)\times\boldsymbol {z_{uw}}
	+\nabla\left(\pi_\rho-\pi\right)=0.$$
Multiplying this equation by $r_\rho$, taking into account (\ref{rm2}), we obtain	
	\begin{align}	\label{energy_r}
	&\tfrac{1}{2}\tfrac{d}{dt}
	 \left(\left\|r_{\rho}(t)\right\|_2^2
	+2\alpha \left\|Dr_{\rho}(t)\right\|_2^2\right)
	+2\nu \left\|D r_\rho(t)\right\|_2^2\nonumber\\
	&=-\left(\boldsymbol{ curl}\,\sigma\left(\boldsymbol r_\rho(t)\right)\times \boldsymbol {y_u}(t),\boldsymbol r_\rho(t)\right)
		-\left( \boldsymbol{ curl}\,\sigma\left(\boldsymbol y_\rho(t)\right)
	\times \boldsymbol r_\rho(t),\boldsymbol r_\rho(t)\right)\nonumber\\
	&-\left(\boldsymbol{ curl}\,\sigma\left(\boldsymbol y_\rho(t)-\boldsymbol {y_u}(t)\right)
	\times \boldsymbol {z_{uw}}(t),\boldsymbol r_\rho(t)\right)\nonumber\\
	&=-\left(\boldsymbol{ curl}\,\sigma\left(\boldsymbol r_\rho(t)\right)\times \boldsymbol {y_u}(t),\boldsymbol r_\rho(t)\right)
	-\left(\boldsymbol{ curl}\,\sigma\left(\boldsymbol y_\rho(t)-\boldsymbol {y_u}(t)\right)
	\times \boldsymbol {z_{uw}}(t),\boldsymbol r_\rho(t)\right),
	\end{align}
where we have used identity (\ref{trilinear_state}). Moreover, due to  (\ref{rm2_lin}), we have
	\begin{align}\left|\left(\boldsymbol{ curl}\,
	\sigma\left(\boldsymbol r_\rho(t)\right)\times  \boldsymbol {y_u}(t),\boldsymbol r_\rho(t)\right)\right|&\leq 
	c\left\|y_u(t)\right\|_{H^3}\left((1+\alpha) \left\|r_\rho(t)\right\|_2^2+\alpha\left\|Dr_\rho(t)\right\|_2^2\right)
		\end{align}
and by using similar arguments, we obtain
	\begin{align}\label{energy_r4}&\left|\left(\boldsymbol{ curl}\,\sigma\left(\boldsymbol y_\rho(t)-\boldsymbol {y_u}(t)\right)
	\times \boldsymbol {z_{uw}}(t),\boldsymbol r_\rho(t)\right)\right|\nonumber\\
	&=\left|b\left(r_\rho(t),z_{uw}(t),\sigma\left(y_\rho(t)-y_u(t)
	\right)\right)-
	b\left(z_{uw}(t),r_\rho(t),\sigma\left(y_\rho(t)-y_u(t)
	\right)\right)\right|
	\nonumber\\
	&\leq \left(\left\|r_\rho(t)\right\|_4 \left\|\nabla z_{uw}(t)
	\right\|_4
+\left\|z_{uw}(t)\right\|_\infty \left\|\nabla r_\rho(t)\right\|_2\right)
	\left\|\sigma\left(y_\rho(t)-y_{u}(t)\right)\right\|_2
	\nonumber\\
	&\leq c\left\|z_{uw}(t)\right\|_{H^2}  \left\|\nabla r_\rho(t)\right\|_2
	\left\|\sigma\left(y_\rho(t)-y_u(t)\right)\right\|_2
	\nonumber\\
	&\leq \tfrac{c}{2}\left\|z_{uw}(t)\right\|^2_{H^2}
	 \left\|\nabla r_\rho(t)\right\|^2_2+
	\tfrac{1}{2}\left\|\sigma\left(y_\rho(t)-y_u(t)\right)\right\|^2_2\nonumber\\
	&\leq c\left\|z_{uw}(t)\right\|^2_{H^2}	\left(\left\|r_\rho(t)\right\|^2_2+ \left\|Dr_\rho(t)\right\|^2_2\right)+
	\tfrac{1}{2}\left\|\sigma\left(y_\rho(t)-y_u(t)\right)\right\|^2_2.
	\end{align}
Combining (\ref{energy_r})-(\ref{energy_r4}) yields
	$$\begin{array}{ll}
	\frac{d}{dt}
	 \left(\left\|r_{\rho}(t)\right\|_2^2
	+2\alpha \left\|Dr_{\rho}(t)\right\|_2^2
	\right)	\vspace{2mm}\\
	\leq \tfrac{c(1+\alpha)^2}{\alpha}
	\left(\left\|y_u(t)\right\|_{H^3}+
	\left\|z_{uw}(t)\right\|^2_{H^2}\right)
		\left(\left\|r_{\rho}(t)\right\|_2^2
	+2\alpha \left\|Dr_{\rho}(t)\right\|_2^2\right)
	+\left\|\sigma\left(y_\rho(t)-y_u(t)\right)\right\|^2_2\vspace{2mm}\\
	\leq \tfrac{c(1+\alpha)^2}{\alpha}
	\left(\left\|y_u\right\|_{L^\infty(I;H^3)}+
	\left\|z_{uw}\right\|^2_{L^\infty(I;H^2)}\right)
		\left(\left\|r_{\rho}(t)\right\|_2^2
	+2\alpha \left\|Dr_{\rho}(t)\right\|_2^2\right)
	+\left\|\sigma\left(y_\rho(t)-y_u(t)\right)\right\|^2_2
	\end{array}$$
and by applying  Gronwall's inequality, we deduce that
	$$\begin{array}{ll}	
	\left\|r_{\rho}\right\|_{L^\infty(I;L^2)}^2
	+2\alpha \left\|Dr_{\rho}\right\|_{L^\infty(I;L^2)}^2\vspace{2mm}\\
	\displaystyle \leq 
	\left\|\sigma\left(y_\rho-y_u\right)\right\|_{2,Q}^2
	{\rm exp}
	\left(\tfrac{cT(1+\alpha)^2}{\alpha}\left(\left\|y_u\right\|_{L^\infty(I;H^3)}+
	\left\|z_{uw}\right\|^2_{L^\infty(I;H^2)}\right)\right).
	\end{array}$$
On the other hand, by taking into account (\ref{yh_sigma}), we have
	$$\begin{array}{ll}\left\|\sigma\left(y_\rho-y\right)\right\|_{2,Q}&\leq \left\|y_\rho-y_u\right\|_{2,Q}+c\alpha 
	\left\|y_\rho-y_u\right\|_{L^2(I;H^2)}\vspace{2mm}\\
	&\leq c\left(\left\|y_\rho-y_u\right\|_{2,Q}+\alpha 
	\left\|\mathbb{A}\left(y_\rho-y_u\right)\right\|_{2,Q}\right)\vspace{2mm}\\
	&\leq cT^\frac{1}{2}\left(\left\|y_\rho-y_u\right\|_{L^\infty(I;L^2)}+\alpha 
	\left\|\mathbb{A}\left(y_\rho-y_u\right)\right\|_{L^\infty(I;L^2)}\right)\end{array}$$
and due to  Proposition \ref{Lipschitz}, we obtain
$$\begin{array}{ll}
\left\| y_{\rho}-y_{u}\right\|_{L^\infty(I;L^2)}^2&\leq \rho^2
\left\| w\right\|^2_{2,Q}
{\rm exp}\left(cT \left(1+\left(1+\alpha \right)\left\|y_{u}\right\|_{L^\infty(I;H^3)}\right)\right)\vspace{2mm}\\
	&\longrightarrow 0 \qquad \mbox{when} \ \rho \rightarrow 0\end{array}$$
and
	$$\alpha 
\left\|\mathbb{A}\left(y_{\rho}-y_{u}\right)
	\right\|_{L^\infty(I;L^2)}^2\leq \left(\tfrac{\rho^2}{\nu}\left\| w\right\|^2_{2,Q}+cT(1+\alpha)M_\rho
\left\|y_{\rho}-y_{u}\right\|^2_{L^\infty(I;L^2)}\right)
	{\rm exp}
\left(\tfrac{c(1+\alpha)T}{\alpha}M_\rho\right),$$
where $c$ is a positive constant independent of $\rho$ and where $M_\rho=\left\|y_{\rho}\right\|_{L^\infty(I;H^3)}+\left\|y_{u}\right\|_{L^\infty(I;H^3)}$. In regards to 
(\ref{state_est012})-(\ref{state_est4}), it follows that $(M_\rho)_\rho$ is uniformly bounded and thus 
	$$\alpha 
\left\|\mathbb{A}\left(y_{\rho}-y_{u}\right)
	\right\|_{L^\infty(I;L^2)}^2\longrightarrow 0 \qquad \mbox{when} \ \rho \rightarrow 0.$$
Combining theses estimates, we deduce that
	$$\lim_{\rho\rightarrow 0}\left(\left\|r_{\rho}\right\|_{L^\infty(I;L^2)}^2
	+2\alpha \left\|Dr_{\rho}\right\|_{L^\infty(I;L^2)}^2\right)=0$$
which gives (\ref{gateau_1}). Identity (\ref{gateau_2}) follows by using standard arguments. $\hfill\Box$

\section{Adjoint equation}\label{sec_adjoint}
\setcounter{equation}{0}
Let $y\in L^\infty(I;W\cap H^3(\Omega))$ and $ f\in L^2(Q)$. The aim of this section is to study the solvability of the adjoint equation defined by 
\begin{equation}\label{adjoint}\left\{
  \begin{array}{ll}
    -\displaystyle \tfrac{\partial  \sigma(\boldsymbol p)}{\partial t}-\nu\Delta \boldsymbol  p-  \boldsymbol{curl}\,\sigma(\boldsymbol y)\times \boldsymbol  p +\boldsymbol{curl}\left(\sigma\left(\boldsymbol y\times \boldsymbol p\right)\right)+\nabla \pi=\boldsymbol f&\quad\mbox{in} \ Q,\vspace{2mm}\\
	{\rm div} \,  \boldsymbol p=0&\quad\mbox{in} \ Q,\vspace{2mm}\\
    \boldsymbol p\cdot  \boldsymbol n=0, \qquad 
	\left(\boldsymbol n\cdot D\boldsymbol p\right)\cdot \boldsymbol \tau=0&\quad\mbox{on}\ \Sigma,
	\vspace{2mm}\\
    \boldsymbol p(T)=0&\quad\mbox{in} \ \Omega.
  \end{array}
\right .\end{equation}
Due to (\ref{trilinear_state}), the term $\boldsymbol{curl}\,\sigma(\boldsymbol y)\times \boldsymbol  p$ can be handled as in the case of the nonlinear term in the state equation. In order to manage the term $\boldsymbol{curl}\left(\sigma\left(\boldsymbol y\times \boldsymbol p\right)\right)$ and to give an adequate variational setting for the adjoint equation, we need the following result.
\begin{lemma} 
Let $y, z$ be in $W\cap H^3(\Omega)$ and $\phi$ be in $W$. Then 
	\begin{equation}\label{trilinear_adjoint}\left( \boldsymbol{curl}\, \sigma\left(\boldsymbol y\times \boldsymbol z\right), \boldsymbol\phi\right)
	=b\left(z,y, \sigma(\phi)\right)-b\left(y,z,\sigma(\phi)\right).\end{equation}
\end{lemma}
{\bf Proof.} As already observed when dealing with the linearized equation, 
some specific difficulties related to the boundary terms arise when considering the Navier-slip conditions. If they are not vanishing, these terms need to be managed and satisfactorily estimated.\vspace{1mm}\\
 We will split the present proof into two steps. 
In order to give a sense to the different boundary terms, we 
 first assume that $y, z$ belong to $W\cap H^4(\Omega)$. We next apply a regularization process to prove that the results are still valid for $y, z\in W\cap H^3(\Omega)$.\vspace{1mm}\\
{\it Step 1.} Standard arguments show that
	\begin{align}\label{curl_sigma}\left( \boldsymbol{curl}\, \sigma\left(\boldsymbol y\times \boldsymbol z\right), \boldsymbol\phi\right)&=
	\left( \boldsymbol{curl}\, \left(\boldsymbol y\times \boldsymbol z\right), \boldsymbol\phi\right)-\alpha
	\left( \boldsymbol{curl}\, \Delta\left(\boldsymbol y\times \boldsymbol z\right), \boldsymbol\phi\right)\nonumber\\
	&=b\left(z,y,\phi\right)-b\left(y,z,\phi\right)-\alpha
	\left( \boldsymbol{curl}\, \Delta\left(\boldsymbol y\times \boldsymbol z\right), \boldsymbol\phi\right),\end{align}
where
	\begin{align}\label{curl_delta}\left(\boldsymbol{curl}\left(\Delta\left(\boldsymbol y\times \boldsymbol z\right)\right), \boldsymbol\phi\right)&=\displaystyle\left(\Delta\left(\boldsymbol y\times \boldsymbol z\right),\boldsymbol{curl}\, \boldsymbol\phi\right)+I_1\nonumber\\
	&=\left(-\boldsymbol{curl}\left(\boldsymbol{curl}\left(\boldsymbol y\times \boldsymbol z\right)\right)
	+\nabla\left( \mathrm{div}\left(\boldsymbol y\times \boldsymbol z\right)\right),\boldsymbol{curl}\, \boldsymbol\phi\right)+I_1\nonumber\\
	&=\left(-\boldsymbol{curl}\left(\boldsymbol{curl}\left(\boldsymbol y\times \boldsymbol z\right)\right),\boldsymbol{curl}\,
	 \boldsymbol\phi\right)+I_1\nonumber\\
	&=-\left(\boldsymbol{curl}\left(\boldsymbol y\times \boldsymbol z\right),\boldsymbol{curl}\left(\boldsymbol{curl}\, \boldsymbol\phi\right)\right)
	\displaystyle -I_2+I_1\nonumber\\
	&=\left(\boldsymbol{curl}\left(\boldsymbol y\times \boldsymbol z\right),\Delta \boldsymbol\phi-\nabla\left( \mathrm{div}\, \boldsymbol\phi\right)\right)+I_1-I_2\nonumber\\
	&=b(\boldsymbol z,\boldsymbol y,\Delta\boldsymbol\phi)-
	b(\boldsymbol y,\boldsymbol z,\Delta \boldsymbol\phi)+I_1-I_2\nonumber\\
	&=b(z,y,\Delta\phi)-b(y,z,\Delta\phi)+I_1-I_2
	\end{align}
with  $$\begin{array}{ll}I_1=\displaystyle \int_\Gamma
	\left(\Delta\left(\boldsymbol y\times \boldsymbol z\right)\times \boldsymbol\phi\right)\cdot \boldsymbol n\,dS
	\qquad\mbox{and}\qquad 
	I_2=\displaystyle \int_\Gamma \boldsymbol{curl}\left(\boldsymbol y\times \boldsymbol z\right)\times \boldsymbol{curl}\,\boldsymbol\phi\cdot \boldsymbol n\,dS. \end{array}$$
Let us now prove that $I_1-I_2=0$.
 By taking into account Lemma \ref{curl_trace}, we have
	\begin{align}\label{I_2}
	\boldsymbol{curl}\left(\boldsymbol y\times \boldsymbol z\right)\times 
	\boldsymbol{curl}\, \boldsymbol \phi \cdot \boldsymbol n\big|_{\Gamma}
	&=\left(\boldsymbol z\cdot \nabla \boldsymbol y-\boldsymbol y\cdot \nabla \boldsymbol z\right)\times 
	(0,0,\phi\cdot g)^\top \cdot \boldsymbol n\big|_{\Gamma}\nonumber\\
	&=\left(\phi\cdot g\right)\left(\boldsymbol z\cdot \nabla \boldsymbol y-\boldsymbol y\cdot \nabla \boldsymbol z\right)\cdot \boldsymbol\tau\big|_{\Gamma}\nonumber\\
	&=\left(\phi\cdot g\right)
	\left(z\cdot \nabla y-y\cdot \nabla z\right)\cdot 
	\tau\big|_{\Gamma}.
	\end{align}
Similarly, since
$$\begin{array}{ll}\Delta\left(\boldsymbol y\times \boldsymbol z\right)&
	=-\boldsymbol{curl}\left(\boldsymbol{curl}\left(\boldsymbol y\times \boldsymbol z\right)\right)
	+\nabla\left( \mathrm{div}\left(\boldsymbol y\times \boldsymbol z\right)\right)=-\boldsymbol{curl}\left(\boldsymbol z\cdot \nabla \boldsymbol y-\boldsymbol y\cdot \nabla \boldsymbol z\right)\vspace{1mm}\\
	&=-\boldsymbol z\cdot \nabla \left(\boldsymbol{curl}\,\boldsymbol y\right)+\boldsymbol{curl}\,\boldsymbol y\cdot \nabla \boldsymbol z-
	\left( \mathrm{div}\,\boldsymbol z\right)\cdot \boldsymbol{curl} \,\boldsymbol y-\displaystyle \sum_{k=1}^2
	\nabla \boldsymbol z_k\times \nabla\boldsymbol y_k\vspace{1mm}\\
	&+\boldsymbol y\cdot \nabla \left(\boldsymbol{curl}\,\boldsymbol z\right)-\boldsymbol{curl}\,\boldsymbol z\cdot \nabla \boldsymbol y+
	\left( \mathrm{div}\,\boldsymbol y\right)\cdot \boldsymbol{curl} \,\boldsymbol z+\displaystyle \sum_{k=1}^2
	\nabla \boldsymbol y_k\times \nabla\boldsymbol z_k\vspace{1mm}\\
	&=\boldsymbol y\cdot \nabla \left(\boldsymbol{curl}\,\boldsymbol z\right)-\boldsymbol z\cdot \nabla \left(\boldsymbol{curl}\,\boldsymbol y\right)+2\displaystyle \sum_{k=1}^2
	\nabla \boldsymbol y_k\times \nabla\boldsymbol z_k\end{array}$$
we deduce that
	$$\begin{array}{ll}\Delta\left(\boldsymbol y\times \boldsymbol z\right)\big|_{\Gamma}-\displaystyle 2\sum_{k=1}^2\left(\nabla \boldsymbol y_k\times \nabla\boldsymbol z_k\right)\big|_{\Gamma}
	&=\boldsymbol y\cdot \nabla \left(\boldsymbol{curl}\,\boldsymbol z\right)\big|_{\Gamma}-\boldsymbol z\cdot \nabla \left(\boldsymbol{curl}\,\boldsymbol y\right)\big|_{\Gamma}\vspace{0mm}\\
	&=\left(0,0,y\cdot \nabla\left(z\cdot g\right)-z\cdot \nabla\left(y\cdot g\right)\right)^\top\big|_{\Gamma}\vspace{2mm}\\
	&=\left(0,0,\left(y\cdot \nabla z\right)\cdot g+\left(y\cdot \nabla g\right)\cdot z\right)^\top\big|_{\Gamma}\vspace{2mm}\\
	&-\left(0,0,-\left(z\cdot \nabla y\right)\cdot g-\left(y\cdot \nabla g\right)\cdot y\right)^\top\big|_{\Gamma}
	\end{array}
	$$
and thus
	\begin{align}\label{J1J2}&\left(\Delta\left(\boldsymbol y\times \boldsymbol z\right)\times \boldsymbol\phi\right)\cdot \boldsymbol n\big|_{\Gamma}-2\displaystyle \sum_{k=1}^2\left(\left(\nabla \boldsymbol y_k\times \nabla\boldsymbol z_k\right)\times \boldsymbol \phi\right)\cdot \boldsymbol n\big|_{\Gamma}\nonumber\\
	&=
	-\left(y\cdot \nabla z-z\cdot \nabla y\right)\cdot g \left(\phi\cdot \tau\right)\big|_{\Gamma}-\left(\left(y\cdot \nabla g\right)
	\cdot z-\left(z\cdot \nabla g\right)\cdot y\right)\left(\phi\cdot \tau\right)\big|_{\Gamma}.\end{align}
Easy calculations, together with the fact that $y\big|_{\Gamma}=
\left(y\cdot \tau\right)\tau\big|_{\Gamma}$ and $z\big|_{\Gamma}=
\left(z\cdot \tau\right)\tau\big|_{\Gamma}$,  show that
	\begin{align}\label{J1}
	\left(y\cdot \nabla g\right)\cdot z\big|_{\Gamma}-
	\left(z\cdot \nabla g\right)\cdot y\big|_{\Gamma}
	&=\left(\nabla g\, y\right)\cdot z\big|_{\Gamma}-
	\left(\nabla g\, z\right)\cdot y\big|_{\Gamma}\nonumber\\
	&=\left(\nabla g\, \left(y\cdot \tau\right)\tau\right)\cdot \left(z\cdot \tau\right)\tau\big|_{\Gamma}-
	\left(\nabla g\, \left(z\cdot \tau\right)\tau\right)\cdot \left(y\cdot \tau\right)\tau\big|_{\Gamma}
	\nonumber\\
	&=\left(y\cdot \tau\right)\left(z\cdot \tau\right)
	\left(\left(\nabla g\, \tau\right)\cdot \tau-
	\left(\nabla g\, \tau\right)\cdot \tau\right)\big|_{\Gamma}=0.
	\end{align}
Taking into account (\ref{I_2})-(\ref{J1}), we deduce that
	$$\begin{array}{ll}
	I_1-I_2&=\displaystyle 2\sum_{k=1}^2\int_\Gamma \left(\nabla \boldsymbol y_k\times \nabla \boldsymbol z_k\right)\times \boldsymbol\phi\cdot \boldsymbol n \,dS\vspace{2mm}\\
	&+\displaystyle\int_\Gamma
	 \left(\left(\phi\cdot \tau\right)\left(z\cdot \nabla y-y\cdot \nabla z\right)\cdot g-\left(\phi\cdot g\right)\left(z\cdot \nabla y-y\cdot \nabla z\right)\cdot \tau\right)dS.
	\end{array}$$
On the other hand, due to Lemma 4.1 in \cite{K06} we have
	$$\left(z\cdot \nabla y-y\cdot \nabla z\right)\cdot n\big |_\Gamma=0.$$
Observing that $u\cdot n\big|_{\Gamma}=0$ implies that
	$$\begin{array}{ll}\left(\phi\cdot \tau\right)\left(u\cdot g\right)
	 -\left(\phi\cdot g\right)\left(u\cdot \tau\right)\big|_{\Gamma}&
	=\left(\phi\cdot \tau\right)\left(\left(u\cdot \tau\right)\tau\cdot g\right)-\left(\left(\phi\cdot \tau\right)\tau\cdot g\right)\left(u\cdot \tau\right)\big|_{\Gamma}
	\vspace{2mm}\\
	&=\left(\phi\cdot \tau\right)\left(u\cdot \tau\right)
	\left(\tau\cdot g-\tau\cdot g\right)\big|_{\Gamma}=0,\end{array}$$
we deduce that
	$$I_1-I_2	=\displaystyle 2\sum_{k=1}^2\int_\Gamma 
	\left(\nabla \boldsymbol y_k\times \nabla \boldsymbol z_k\right)
	\times \boldsymbol\phi\cdot \boldsymbol n \,dS.$$
Finally, since $\left(Du\cdot n\right)\cdot \tau\big|_{\Gamma}=0$ and $\mathrm{div}\, u=0$ imply 
	$$\left(n_2^2-n_1^2\right)\left(\tfrac{\partial u_1}{\partial x_2}+
	\tfrac{\partial u_2}{\partial x_1}\right)-2n_1n_2\,\tfrac{\partial u_1}{\partial x_1}=0 \qquad \mbox{on} \ \Gamma,$$
we obtain
	$$\begin{array}{ll}
	\displaystyle \sum_{k=1}^2\left(\nabla \boldsymbol y_k\times \nabla\boldsymbol z_k\right)\times \boldsymbol \phi\cdot \boldsymbol n\big|_{\Gamma}
	&\displaystyle =\sum_{k=1}^2\left(\tfrac{\partial y_k}{\partial x_1}\tfrac{\partial z_k}{\partial x_2}-
	\tfrac{\partial y_k}{\partial x_2}\tfrac{\partial z_k}{\partial x_1}\right) \phi\cdot \tau\big|_{\Gamma}\vspace{2mm}\\
	&\displaystyle =\left(\tfrac{\partial y_1}{\partial x_1}\left(\tfrac{\partial z_1}{\partial x_2}+
	\tfrac{\partial z_2}{\partial x_1}\right)-\tfrac{\partial z_1}{\partial x_1}\left(\tfrac{\partial y_1}{\partial x_2}+
	\tfrac{\partial y_2}{\partial x_1}\right)\right) \phi\cdot \tau\big|_{\Gamma}=0
\end{array}$$
and thus
	\begin{equation}\label{I1-I2}I_1-I_2=0.\end{equation}
The conclusion is then a consequence of (\ref{curl_sigma}), (\ref{curl_delta}) and (\ref{I1-I2}).\vspace{1mm}\\
{\it Step 2. Regularization process.} Let us now go back to the case 
$y, z\in W\cap H^3(\Omega)$. We first infer that there exist $y_\varepsilon, z_\varepsilon \in W\cap H^4(\Omega)$ such that
	$$\lim_{\varepsilon\rightarrow 0^+}
	\left\|y-y_\varepsilon\right\|_{H^3}
	=\lim_{\varepsilon\rightarrow 0^+}
	\left\|z-z_\varepsilon\right\|_{H^3}=0.$$
Indeed, if $y\in W\cap H^3(\Omega)$, then it  satisfies the following Stokes system
	\begin{equation}\label{stokes_molifier}\left\{\begin{array}{ll}
   	 -\Delta y+y+\nabla \pi=f&\quad\mbox{in} \ \Omega,\vspace{2mm}\\
	 \mathrm{div} \,  y=0&\quad\mbox{in} \ \Omega,\vspace{2mm}\\
   	 y\cdot n=0, \qquad 
	\left(n\cdot D y\right)\cdot \tau=0&\quad\mbox{on}\ \Gamma
 	 \end{array}
	\right. \end{equation}
with $f=-\Delta y+y\in H^1(\Omega)$ and $\pi=0$. Using Friedrichs mollifiers, we can construct 
$f_\varepsilon\in H^2(\Omega)$ such that
	$$\lim_{\varepsilon \rightarrow 0^+}
	\left\|f-f_\varepsilon\right\|_{H^1}=0.$$
Let $y_\varepsilon\in H^4(\Omega)$ be the solution of (\ref{stokes_molifier}) corresponding to $f_\varepsilon$. By using classical regularity results, we have
	$$\left\|y-y_\varepsilon\right\|_{H^3} \leq c 
	\left\|f-f_\varepsilon\right\|_{H^1} \longrightarrow 0 \quad \mbox{when} \ \varepsilon \rightarrow 0.$$
On the other hand, by taking into account the first step, we deduce that
	$$\left( \boldsymbol{curl}\, \sigma\left(\boldsymbol y_\varepsilon\times \boldsymbol z_\varepsilon\right), \boldsymbol\phi\right)
	=b\left(z_\varepsilon,y_\varepsilon, \sigma(\phi)\right)-b\left(y_\varepsilon,z_\varepsilon,\sigma(\phi)\right)$$
and the result follows by passing to the limit.$\hfill\Box$\vspace{1mm}\\
The identities (\ref{trilinear_state}) and (\ref{trilinear_adjoint}) motivate the following variational formulation for the adjoint equation.
\begin{definition}  A function 
$p\in L^\infty(I;W)$ with  $\frac{\partial p}{\partial t}\in L^2(I;V)$  is a solution of $(\ref{adjoint})$ if
$p(T)=0$ and 
\begin{align}
\label{form_var_lin_adj}
-\left(\tfrac{\partial p(t) }{\partial t}, \phi\right)&
-2\alpha\left(D\tfrac{\partial p(t)}{\partial t},
D\phi\right)+
2\nu\left(Dp(t),D  \phi\right)+
b\left(p(t),\phi, \sigma(y(t)) \right)-b\left( \phi,p(t), \sigma(y(t)) \right)
\vspace{2mm}\nonumber\\
&+
b\left(p(t), y(t), \sigma(\phi) \right)-b\left( y(t),p(t), \sigma(\phi) \right)=
\left(f(t), \phi\right) \qquad \mbox{for all} \  \phi\in  W.
\end{align}
\end{definition}
The solution of $(\ref{adjoint})$ in the sense of this definition is constructed using the Galerkin's method defined in the previous section. Existence of an approximate solution and a corresponding a priori $H^1$ estimate can be established using standard arguments. In order to establish the corresponding $H^2$ in space estimate (and consequently the time derivative estimates neccessary to pass to the limit), we need the following result. 
\begin{lemma}\label{prop_a_adj_3}
Let $y,\psi \in  W\cap H^3(\Omega)$ and $\phi\in  V$. Then
	$$\left|\left(\boldsymbol{curl}\left(\sigma\left(\boldsymbol  y\times \boldsymbol  \psi\right)\right)-\boldsymbol{curl}\,\sigma\left(\boldsymbol  y\right)\times \boldsymbol  \psi,\boldsymbol\phi\right)\right|\leq c(1+\alpha) \left\|y\right\|_{H^3}\left\|\psi\right\|_{H^2}
	\left\|\phi\right\|_{2}
	+\alpha \left|b\left(y,\Delta \psi,\phi\right)\right|,$$
where $c$ is a positive constant depending only on $\Omega$. 
\end{lemma}
{\bf Proof.} Standard calculations show that
	$$\Delta\left(\boldsymbol  y\times \boldsymbol  \psi\right)=\boldsymbol  y\times \Delta \boldsymbol  \psi-\boldsymbol  \psi\times \Delta\boldsymbol  y-
	2\sum_{i=1}^2\tfrac{\partial \boldsymbol  \psi}{\partial x_i}\times \tfrac{\partial \boldsymbol  y}{\partial x_i}$$
and thus 
	$$\begin{array}{ll}\boldsymbol{curl}\left(\Delta\left(\boldsymbol  y\times \boldsymbol  \psi\right)\right)&=\displaystyle
	\boldsymbol{curl}\left(\boldsymbol  y\times \Delta\boldsymbol  \psi\right)-
	\boldsymbol{curl}\left(\boldsymbol  \psi\times \Delta\boldsymbol  y\right)-2 \displaystyle
	\sum_{i=1}^2\boldsymbol{curl}\left(\tfrac{\partial \boldsymbol  \psi}{\partial x_i}\times \tfrac{\partial \boldsymbol  y}{\partial x_i}\right)
	\vspace{1mm}\\
	&=\displaystyle \left( \mathrm{div}\,\Delta\boldsymbol  \psi\right) \boldsymbol  y+\Delta\boldsymbol  \psi\cdot \nabla \boldsymbol  y-\left( \mathrm{div}\,\boldsymbol  y\right) \Delta\boldsymbol  \psi-\boldsymbol  y\cdot \nabla \left(\Delta\boldsymbol  \psi\right)\vspace{2mm}\\
	&-\displaystyle \left( \mathrm{div}\,\Delta\boldsymbol  y\right) \boldsymbol  \psi-\Delta\boldsymbol  y\cdot \nabla \boldsymbol  \psi+\left( \mathrm{div}\,\boldsymbol  \psi\right) \Delta\boldsymbol  y+\boldsymbol  \psi\cdot \nabla \left(\Delta\boldsymbol  y\right)\vspace{2mm}\\
	&\displaystyle-2\displaystyle\sum_{i=1}^2\left( \mathrm{div}\left(\tfrac{\partial \boldsymbol  \psi}{\partial x_i}\right)
	\tfrac{\partial \boldsymbol  y}{\partial x_i}+\tfrac{\partial \boldsymbol  \psi}{\partial x_i}\cdot \nabla \left(\tfrac{\partial \boldsymbol  y}{\partial x_i}\right)
	- \mathrm{div}\left(\tfrac{\partial \boldsymbol  y}{\partial x_i}\right)
	\tfrac{\partial \boldsymbol  \psi}{\partial x_i}-\tfrac{\partial \boldsymbol  y}{\partial x_i}\cdot \nabla \left(\tfrac{\partial \boldsymbol  \psi}{\partial x_i}\right)\right)
	\vspace{2mm}\\
	&=\displaystyle \Delta\boldsymbol  \psi\cdot \nabla \boldsymbol  y-\boldsymbol  y\cdot \nabla \left(\Delta\boldsymbol  \psi\right)-\Delta\boldsymbol  y\cdot \nabla \boldsymbol  \psi+\boldsymbol  \psi\cdot \nabla \Delta\boldsymbol  y\vspace{2mm}\\
	&\displaystyle-2\displaystyle\sum_{i=1}^2\left(
	\tfrac{\partial \boldsymbol  \psi}{\partial x_i}\cdot \nabla \left(\tfrac{\partial \boldsymbol  y}{\partial x_i}\right)
	-\tfrac{\partial \boldsymbol  y}{\partial x_i}\cdot \nabla \left(\tfrac{\partial \boldsymbol  \psi}{\partial x_i}\right)\right).\end{array}$$
Therefore
	$$\begin{array}{ll}\left(\boldsymbol{curl}\left(\Delta\left(\boldsymbol  y\times \boldsymbol  \psi\right)\right),\boldsymbol\phi\right)&=
	b\left(\Delta \psi,y,\phi\right)-b\left(y,\Delta \psi,\phi\right)-b\left(\Delta y,\psi,\phi\right)
	+b\left(\psi,\Delta y,\phi\right)\vspace{2mm}\\
	&-2\displaystyle\sum_{i=1}^2\left(b\left(
	\tfrac{\partial \psi}{\partial x_i},\tfrac{\partial y}{\partial x_i},\phi\right)
	-b\left(\tfrac{\partial y}{\partial x_i},\tfrac{\partial \psi}{\partial x_i},\phi\right)\right).
	\end{array}$$
On the other hand, by taking into account (\ref{trilinear_state}), we have
	$$\left(\boldsymbol{curl}\left(\Delta\boldsymbol  y\right)\times \boldsymbol  \psi,\boldsymbol\phi\right)=b(\phi,\psi,\Delta y)-b(\psi,\phi,\Delta y)=b(\phi,\psi,\Delta y)+b(\psi,\Delta y,\phi).$$
Combining the previous two identities, we obtain
	$$\begin{array}{ll}\left(\boldsymbol{curl}\left(\Delta\left(\boldsymbol  y\times \boldsymbol  \psi\right)\right)-\boldsymbol{curl}\left(\Delta\boldsymbol  y\right)\times \boldsymbol  \psi,\boldsymbol\phi\right)
	&=b\left(\Delta \psi,y,\phi\right)-b\left(y,\Delta \psi,\phi\right)-b\left(\Delta y,\psi,\phi\right)-b(\phi,\psi,\Delta y)\vspace{2mm}\\
	&-2\displaystyle\sum_{i=1}^2\left(b\left(
	\tfrac{\partial \psi}{\partial x_i},\tfrac{\partial y}{\partial x_i},\phi\right)
	-b\left(\tfrac{\partial y}{\partial x_i},\tfrac{\partial \psi}{\partial x_i},\phi\right)\right),\end{array}$$
which together with the fact that
	$$\begin{array}{ll}\left(\boldsymbol{curl}\left(\boldsymbol  y\times \boldsymbol  \psi\right),\boldsymbol\phi\right)-\left(\boldsymbol{curl}\,\boldsymbol  y\times \boldsymbol  \psi,\boldsymbol\phi\right)&=\left(\boldsymbol{curl}\left(\boldsymbol  y\times \boldsymbol  \psi\right),\boldsymbol\phi\right)-\left(\boldsymbol  y,
	\boldsymbol{curl}\left(\boldsymbol  \psi\times\boldsymbol\phi\right)\right)\vspace{2mm}\\
	&=\left(\boldsymbol  \psi\cdot \nabla \boldsymbol  y-
	\boldsymbol  y\cdot \nabla \boldsymbol  \psi,\boldsymbol\phi\right)-\left(\boldsymbol  y,\boldsymbol\phi\cdot \nabla \boldsymbol  \psi-\boldsymbol  \psi\cdot \nabla \boldsymbol\phi
	\right)\vspace{2mm}\\
	&=b\left(\psi,y,\phi\right)-b\left(y,\psi,\phi\right)-b\left(\phi,\psi,y\right)+
	b\left(\psi,\phi,y\right)\vspace{2mm}\\
	&=-b\left(y,\psi,\phi\right)-	b\left(\phi,\psi,y\right)\vspace{2mm}\\
	&=-b\left(y,\psi,\phi\right)+b\left(\phi,y,\psi\right)\end{array}$$
gives 
$$\begin{array}{ll}
	\left(\boldsymbol{curl}\left(\sigma\left(\boldsymbol  y\times \boldsymbol  \psi\right)\right)-\boldsymbol{curl}\,\sigma\left(\boldsymbol  y\right)\times \boldsymbol  \psi,\boldsymbol\phi\right)
	&=b\left(\phi,y,\psi\right)-b\left(y,\psi,\phi\right)-\alpha b\left(\Delta \psi,y,\phi\right)+\alpha b\left(y,\Delta \psi,\phi\right)\vspace{2mm}\\
	&+\alpha b\left(\Delta y,\psi,\phi\right)+\alpha b(\phi,\psi,\Delta y)\vspace{0mm}\\
	&+2\alpha\displaystyle\sum_{i=1}^2\left(b\left(
	\tfrac{\partial \psi}{\partial x_i},\tfrac{\partial y}{\partial x_i},\phi\right)
	-b\left(\tfrac{\partial y}{\partial x_i},\tfrac{\partial \psi}{\partial x_i},\phi\right)\right).
	\end{array}$$
Therefore,
	$$\left|\left(\boldsymbol{curl}\left(\sigma\left(\boldsymbol  y\times \boldsymbol  \psi\right)\right)-\boldsymbol{curl}\,\sigma\left(\boldsymbol  y\right)\times \boldsymbol  \psi,\boldsymbol\phi\right)\right|
	$$
	$$\leq \left|b\left(\phi,y,\psi\right)\right|+\left|b\left(y,\psi,\phi\right)\right|+\alpha \left|b\left(\Delta \psi,y,\phi\right)\right|
	+\alpha \left|b\left(\Delta y,\psi,\phi\right)\right|+\alpha \left|b(\phi,\psi,\Delta y)\right|$$
	$$+\alpha \left|b\left(y,\Delta \psi,\phi\right)\right|+2\alpha\displaystyle\sum_{i=1}^2\left(\left|b\left(
	\tfrac{\partial \psi}{\partial x_i},\tfrac{\partial y}{\partial x_i},\phi\right)\right|
	+\left|b\left(\tfrac{\partial y}{\partial x_i},\tfrac{\partial \psi}{\partial x_i},\phi\right)\right|\right)$$
	$$\leq \left(\left\|\nabla y\right\|_4\left\|\psi\right\|_4+\left\|y\right\|_4
	\left\|\nabla \psi\right\|_4+\alpha \left\|\nabla y\right\|_\infty\left\|\Delta \psi\right\|_2+2\alpha \left\|\Delta y\right\|_4
	\left\|\nabla \psi\right\|_4\right)
	\left\|\phi\right\|_2$$
	$$+\alpha \left|b\left(y,\Delta \psi,\phi\right)\right|+c\alpha\left(
	\left\|\nabla y\right\|_\infty\left\|\psi\right\|_{H^2}+\left\|y\right\|_{2,4}\left\|\nabla \psi\right\|_4\right)
	\left\|\phi\right\|_2$$
	$$\leq c(1+\alpha)\left\|y\right\|_{H^3}
	\left\|\psi\right\|_{H^2}\left\|\phi\right\|_2
	+\alpha \left|b\left(y,\Delta \psi,\phi\right)\right|$$
which gives the claimed result.$\hfill\Box$\vspace{1mm}\\
The next result deals with the existence and uniqueness of a solution for the adjoint equation.
\begin{proposition} \label{ex_uniq_adj}  Assume that $ f\in L^2(Q)$. Then equation $(\ref{adjoint})$ admits a unique solution
 $p\in L^\infty(I;W)$ with  $\frac{\partial p}{\partial t}\in L^2(I;V)$. 
  Moreover, the following estimates hold
$$\left\|p\right\|_{L^\infty(I;L^2)}^2+2\alpha
\left\|Dp\right\|_{L^\infty(I;L^2)}^2\leq
\left\|f\right\|^2_{2,Q}
\displaystyle
e^{cT \left(1+\left(1+\alpha \right)\left\|y\right\|_{L^\infty(I;H^3)}\right)},$$		
$$\left\|Dp\right\|_{L^\infty(I;L^2)}^2+\alpha 
\left\|\mathbb{A}p\right\|_{L^\infty(I;L^2)}^2$$
$$\leq \left(\tfrac{1}{\nu}\left\|f\right\|^2_{2,Q}+cT(1+\alpha)\left\|y\right\|_{L^\infty(I;H^3)}
\left\|p\right\|^2_{L^\infty(I;L^2)}\right)e^{
\frac{c(1+\alpha)T}{\alpha}\left\|y\right\|_{L^\infty(I;H^3)}},$$	
$$\left\|\tfrac{\partial p}{\partial t}\right\|_{2,Q}^2+
\alpha\left\|D\tfrac{\partial p}{\partial t}\right\|_{2,Q}^2\leq 
c\left(\left\|f\right\|_{2,Q}^2+T\left(\nu^2+\tfrac{(1+\alpha)^2}{\alpha}\left\|y\right\|_{L^\infty(I;H^3)}^2\right)
	\left\|p\right\|_{L^\infty(I;H^2)}^2\right),$$
where $c$ is a constant only depending on $\Omega$. 
\end{proposition}
{\bf Proof.} Let us first notice that $p$ is  solution of the terminal value problem (\ref{adjoint}) if and only if 
$\psi$ defined as $\psi(t)=p(T-t)$ is the solution of the following initial value problem
\begin{equation}\label{adjointF}\left\{
  \begin{array}{ll}
    \displaystyle \tfrac{\partial  \sigma(\boldsymbol\psi)}{\partial t}-\nu\Delta \boldsymbol  \psi-  \boldsymbol{curl}\,\sigma(\tilde{\boldsymbol y})\times \boldsymbol\psi +\boldsymbol{curl}\left(\sigma\left(\tilde{\boldsymbol y}\times \boldsymbol\psi\right)\right)+\nabla\tilde{ \pi}=\tilde{\boldsymbol f}&\quad\mbox{in} \ Q,\vspace{2mm}\\
	{\rm div} \, \boldsymbol \psi=0&\quad\mbox{in} \ Q,\vspace{2mm}\\
   \boldsymbol \psi\cdot \boldsymbol  n=0, \qquad 
	\left(\boldsymbol n \cdot D\boldsymbol\psi \right)\cdot \boldsymbol\tau=0&\quad\mbox{on}\ \Sigma,
	\vspace{2mm}\\
  \boldsymbol\psi(0)=0&\quad\mbox{in} \ \Omega,
  \end{array}
\right .\end{equation}
where $\tilde{y}(t)=y(T-t)$, $\tilde{f}(t)=f(T-t)$ and $\tilde{\pi}(t)=\pi(T-t).$ In regards to (\ref{form_var_lin_adj}), 
$\psi\in L^2(I;W)$ with  $\frac{\partial \psi}{\partial t}\in L^2(I;V)$  is a solution of $(\ref{adjointF})$ if
$\psi(0)=0$ and 
\begin{align}
\label{form_var_lin_adj_psi}
\left(\tfrac{\partial \psi(t) }{\partial t}, \phi\right)&
+2\alpha\left(D\tfrac{\partial\psi(t)}{\partial t},
D\phi\right)+
2\nu\left(D \psi(t),D  \phi\right)+
b\left(\psi(t),\phi, \sigma(\tilde y(t)) \right)-b\left( \phi,\psi(t), \sigma(\tilde y(t)) \right)
\vspace{2mm}\nonumber\\
&+
b\left(\psi(t), \tilde y(t), \sigma(\phi) \right)-b\left( \tilde y(t),\psi(t), \sigma(\phi) \right)=
\left( \tilde f(t), \phi\right) \qquad \mbox{for all} \  \phi\in  W.
\end{align}
The rest of the proof is devoted to the solvability of (\ref{adjointF}). The corresponding approximate problem reads as
\begin{equation}\label{faedo_galerkin_adj}
\left\{\begin{array}{llll}\displaystyle \mbox{Find} \ \psi_k(t)=\sum_{i=1}^k \zeta_{i}(t)
e_i \ \mbox{solution}, \ \mbox{for} \ 
1\leq j\leq k, \ \mbox{of} \vspace{2mm}\\ 
\displaystyle	\left(\tfrac{\partial \psi_k(t) }{\partial t}, e_j\right)+2\alpha
\displaystyle\left(D\tfrac{\partial\psi_k(t) }{\partial t},
De_j\right)+2\nu\left(D\psi_k(t),De_j\right)\vspace{2mm}\\
\ \quad \qquad \qquad+
b\left(\psi_k(t),e_j, \sigma(\tilde y(t)) \right)-b\left(e_j,\psi_k(t), \sigma(\tilde y(t)) \right)\vspace{2mm}\\
\ \quad \qquad \qquad+
b\left(\psi_k(t),\tilde  y(t), \sigma(e_j) \right)-b\left(\tilde y(t),\psi_k(t), \sigma(e_j) \right)
=\left(\tilde f(t), e_j\right),\vspace{2mm}\\
\psi_k(0)=0.
\end{array}\right.
\end{equation}
where $(e_j)_j$ is the basis defined by (\ref{eigen_functions}). This is a linear differential system for $\zeta (t)=(\zeta_1 (t),\dots, \zeta_k (t) )^\top$, of the form
	$$
	A\tfrac{d\zeta(t)}{dt} +\tilde N(t)\zeta (t)=\tilde F(t),
	$$
where $A$ is the nonsingular constant matrix defined by
	$$A_{ij}=\left(e_i,e_j\right)+2\alpha\left(De_i,De_j\right) \qquad i,j=1,\cdots,k,$$
 $\tilde N(t)$ is the matrix given by
	  $$\begin{array}{ll}\tilde N_{ij}(t)&=2\nu\left(De_i,De_j\right)-b\left(e_i,e_j,\sigma(\tilde y(t))\right)+
	b\left(e_j,e_i,\sigma(\tilde y(t))\right)\vspace{2mm}\\
	&+b\left(e_j,\tilde y(t),\sigma(e_i)\right)-b\left(\tilde y(t),e_j,\sigma(e_i)\right) \qquad i,j=1,\cdots,k\end{array}$$ 
and 
	$$\tilde F_i(t)=\left(\tilde f(t),e_i\right)\qquad i=1,\cdots,k.$$
 Using the same arguments as for the linearized equation, we can justify that this system has a unique  solution $\zeta \in H^1(0,T).$ \vspace{2mm}\\
The remaing proof is split into four steps.\vspace{2mm}\\
{\it Step 1. Uniform  $H^1$ in space estimate}.
 Multiplying the equation (\ref{faedo_galerkin_adj}) by $\zeta_j (t)$ and taking the sum over $j=1,\dots,k$, we verify that $\psi_k(t)$ is 
a solution of the following equation 
$$\begin{array}{ll}
	&\displaystyle\tfrac{1}{2}\,\tfrac{d }{d t}\left(\left\| \psi_k(t)\right\|_2^2+2\alpha
\left\|D\psi_k(t)\right\|_2^2\right)+2\nu \left\|D \psi_k(t)\right\|_2^2
	\vspace{2mm}\\
	&=\left(\tilde f(t), \psi_k(t)\right)+
	\left(\boldsymbol{ curl}\,\sigma(\tilde {\boldsymbol y}(t))\times  \boldsymbol \psi_k(t), \boldsymbol\psi_k(t)\right)-
	 \left(\boldsymbol{curl}\,\sigma(\tilde {\boldsymbol y}(t)\times
	  \boldsymbol \psi_k(t)),\boldsymbol\psi_k(t)\right).
	\end{array}$$
Due to (\ref{trilinear_state}) and (\ref{trilinear_adjoint}),
 we have
$$\begin{array}{ll}
	&\left( \boldsymbol{curl}\,\sigma(\tilde {\boldsymbol y}(t))\times\boldsymbol\psi_k(t),\boldsymbol\psi_k(t)\right)
	-\left(\boldsymbol{curl}\left(\sigma\left(\tilde {\boldsymbol y}(t)\times  \boldsymbol\psi_k(t)\right)\right),\boldsymbol\psi_k(t)\right)\vspace{3mm}\\
	&=
	-b\left( \psi_k(t), \tilde y(t), \sigma(\psi_k(t)) \right)+b\left( \tilde y(t), \psi_k(t), \sigma(\psi_k(t)) \right)=-
	\left( \boldsymbol{curl}\,\sigma(\boldsymbol\psi_k(t))\times  \tilde {\boldsymbol y}(t),\boldsymbol\psi_k(t)\right)
\end{array}$$
and thus
$$\begin{array}{ll}
\displaystyle\tfrac{1}{2}\,\tfrac{d }{d t}\left(\left\|  \psi_k(t)\right\|_2^2+2\alpha
\left\|D\psi_k(t)\right\|_2^2\right)+2\nu \left\|D  \psi_k(t)\right\|_2^2
\vspace{2mm}\\
=\left( \tilde f(t),\psi_k(t)\right)-\left( \boldsymbol{ curl}\,\sigma(\boldsymbol \psi_k(t))\times \tilde {\boldsymbol y}(t),
\boldsymbol\psi_k(t)\right).
\end{array}$$ 
Arguing exactly as in the first step of the proof of  Proposition \ref{ex_uniq_lin}, we obtain
\begin{equation}\label{est_psik_H1}
\left\|  \psi_k\right\|_{L^\infty(I;L^2)}^2+2\alpha
\left\|D\psi_k\right\|_{L^\infty(I;L^2)}^2\leq
\left\|\tilde F_1\right\|_{L^1(0,T)}
\displaystyle e^{\tilde G_1T},
\end{equation}
where $\tilde F_1(t)=\left\| \tilde f(t)\right\|^2_2$ and $\tilde G_1=c \left(1+(1+\alpha)
\|\tilde y\|_{L^\infty(I;H^3)}\right)$.\vspace{2mm}\\
 {\it Step 2.  Uniform $H^2$ in space estimate.} Arguing as in the second step of the proof of Proposition \ref{ex_uniq_lin}, we deduce that
$$\begin{array}{ll}
&\displaystyle\tfrac{1}{2}\tfrac{d}{d t}\left(2\left\|D\psi_k(t)\right\|_2^2+\alpha \left\|\mathbb{A}\psi_k(t)\right\|_2^2\right)
+\nu \left\|\mathbb{A}\psi_k(t)\right\|_2^2\vspace{2mm}\\
	&=\left(\tilde f(t),\mathbb{A}\psi_k(t)\right)+\left(
	\boldsymbol{curl}\, \sigma(\tilde {\boldsymbol y}(t))\times
	  \boldsymbol \psi_k(t)-\boldsymbol{curl}\left(\sigma\left(\tilde {\boldsymbol y}(t)\times  \boldsymbol\psi_k(t)\right)\right),\mathbb{A}\boldsymbol\psi_k(t)\right).
	\end{array}$$
Using the Young inequality and Lemma \ref{prop_a_adj_3}, we obtain
$$\begin{array}{ll}
\tfrac{d}{d t}\left(2\left\|D\psi_k(t)\right\|_2^2+\alpha \left\|\mathbb{A}\psi_k(t)\right\|_2^2\right)
+\nu \left\|\mathbb{A}\psi_k(t)\right\|_2^2\vspace{2mm}\\
\leq \tfrac{1}{\nu}\left\|\tilde  f(t)\right\|^2_2+2\left|\left(
	\boldsymbol{curl}\, \sigma(\tilde {\boldsymbol y}(t))\times
	  \boldsymbol \psi_k(t)-\boldsymbol{curl}\left(\sigma\left(\tilde {\boldsymbol y}(t)\times  \boldsymbol\psi_k(t)\right)\right),\mathbb{A}\boldsymbol\psi_k(t)\right)\right|\vspace{2mm}\\
\leq \tfrac{1}{\nu}\left\|\tilde  f(t)\right\|^2_2
	+c(1+\alpha) \left\|\tilde y(t)\right\|_{H^3}
	\left\|\psi_k(t)\right\|_{H^2}
	\left\|\mathbb{A}\boldsymbol\psi_k(t)\right\|_{2}
	+\alpha \left|b\left(\tilde y(t),\Delta \psi_k(t),\mathbb{A}\boldsymbol\psi_k(t)\right)\right|.\end{array}$$
Since
	$$b\left(\tilde y(t),\Delta \psi_k(t),\mathbb{A}\boldsymbol\psi_k(t)\right)=b\left(\tilde y(t),\Delta \psi_k(t)+\mathbb{A}\boldsymbol\psi_k(t),\mathbb{A}\boldsymbol\psi_k(t)\right)$$
by taking into account (\ref{sigma_psigma2}) and (\ref{yh_sigma}), we deduce that
$$\begin{array}{ll}
\tfrac{d}{d t}\left(2\left\|D\psi_k(t)\right\|_2^2+\alpha \left\|\mathbb{A}\psi_k(t)\right\|_2^2\right)
+\nu \left\|\mathbb{A}\psi_k(t)\right\|_2^2\vspace{2mm}\\
\leq \tfrac{1}{\nu}\left\|\tilde  f(t)\right\|^2_2
	+c(1+\alpha) \left\|\tilde y(t)\right\|_{H^3}
	\left\|\psi_k(t)\right\|_{H^2}
	\left\|\mathbb{A}\boldsymbol\psi_k(t)\right\|_{2}
\vspace{2mm}\\
\leq \tfrac{1}{\nu}\left\|\tilde  f(t)\right\|^2_2+c(1+\alpha)\left\|\tilde y\right\|_{H^3}\left(\left\|\psi_k\right\|^2_{2}+\left\|\mathbb{A}\psi_k(t)\right\|^2_{2}
\right) \vspace{2mm}\\
\leq \tilde F_2(t)+\tilde G_2
\left(2\left\|D\psi_k(t)\right\|_2^2+\alpha 
\left\|\mathbb{A}\psi_k(t)\right\|_2^2\right),
	\end{array}$$
where
\begin{align}
&\tilde F_2(t)= \tfrac{1}{\nu}\left\|\tilde f(t)\right\|_2^2+c(1+\alpha)\left\|\tilde y\right\|_{L^\infty(I;H^3)}
\left\|\psi_k\right\|^2_{L^\infty(I;L^2)},\nonumber\\
&\tilde G_2=\tfrac{c(1+\alpha)}{\alpha}\left\|\tilde y\right\|_{L^\infty(I;H^3)}.	\nonumber
		\end{align}
Upon integration, we have
	$$2\left\|D\psi_k(t)\right\|_2^2+\alpha 
\left\|\mathbb{A}\psi_k(t)\right\|_2^2\leq 
	\left\|\tilde F_2\right\|_{L^1(0,t)}+\tilde G_2\int_0^t 
\left(2\left\|D\psi_k(s)\right\|_2^2+\alpha 
\left\|\mathbb{A}\psi_k(s)\right\|_2^2\right)\,ds$$
and by using Gronwall's inequality, we deduce that
\begin{align}
\label{lipschitz_est_V2_psim}
2\left\|D\psi_k\right\|_{L^\infty(I;L^2)}^2+\alpha 
\left\|\mathbb{A}\psi_k\right\|_{L^\infty(I;L^2)}^2
&\leq \left\|\tilde F_2\right\|_{L^1(0,T)} e^{\tilde G_2T}.
\end{align}	
{\it Step 3. Uniform estimate for the  time derivative.} Let us multiply both sides of (\ref{faedo_galerkin_adj}) 
	by  $\frac{d\zeta_{j}(t)}{dt}$ and sum over $j=1, \dots,k$. This gives
$$\begin{array}{ll}
\left\|\tfrac{\partial \psi_k(t) }{\partial t}\right\|_2^2+
\alpha \left\|D\tfrac{\partial \psi_k(t) }{\partial t}\right\|_2^2&=\left(\tilde f(t)+\nu \Delta \psi_k(t), 
\tfrac{\partial \psi_k(t) }{\partial t}\right)\vspace{2mm}\\
&+\left(\boldsymbol{ curl}\,\sigma(\tilde {\boldsymbol y}(t))\times
\boldsymbol \psi_k(t)- \boldsymbol{curl}\left(\sigma\left(\tilde {\boldsymbol y}(t)\times\boldsymbol\psi_k(t)\right)\right),\tfrac{\partial \boldsymbol \psi_k(t) }{\partial t} \right).
\end{array}$$
Due to inequality (\ref{prop_a_adj_3}), we have
	$$\begin{array}{ll}
	&\left|	\left(\boldsymbol{curl}\,\sigma(\tilde {\boldsymbol y}(t))\times
	  \boldsymbol \psi_k(t)-\boldsymbol{curl}\left(\sigma\left(\tilde {\boldsymbol y}(t)\times\boldsymbol\psi_k(t)\right)\right),
\tfrac{\partial \boldsymbol \psi_k(t) }{\partial t}\right)\right|
	\vspace{2mm}\\
	&\leq c(1+\alpha)\left\|\tilde y(t)\right\|_{H^3}\left\|\psi_k(t)\right\|_{H^2}\left\|\tfrac{\partial \psi_k(t)}{\partial t}\right\|_{2}+
 \alpha\left|-b\left(\tilde y(t),\tfrac{\partial \psi_k(t)}{\partial t},\Delta \psi_k(t)\right)   \right|\vspace{2mm}\\
	  &\leq c(1+\alpha)\left\|\tilde y(t)\right\|_{H^3}\left\|\psi_k(t)\right\|_{H^2}\left\|\tfrac{\partial \psi_k(t)}{\partial t}\right\|_{H^1}\vspace{2mm}\\
	  &\leq c(1+\alpha)\left\|\tilde y(t)\right\|_{H^3}\left\|\psi_k(t)\right\|_{H^2}\left(\left\|\tfrac{\partial \psi_k(t)}{\partial t}\right\|_{2}+\left\|D\tfrac{\partial \psi_k(t)}{\partial t}\right\|_{2}\right).
	\end{array}$$
Therefore, by arguing as in the third step of the proof of Proposition \ref{ex_uniq_lin}, we deduce that
\begin{align}
&\left\|\tfrac{\partial  \psi_k}{\partial t}\right\|_{2,Q}^2+
\alpha\left\|D\tfrac{\partial  \psi_k}{\partial t}\right\|_{2,Q}^2\nonumber\\
&\leq 
c\left(\|\tilde f\|_{2,Q}^2+T\left(\nu^2+\tfrac{(1+\alpha)^2}{\alpha}\left\|\tilde y\right\|_{L^\infty(I;H^3)}^2\right)
	\left\| \psi_k\right\|_{L^\infty(I;H^2)}^2\right).
\label{lipschitz_est_V2_psimt}
\end{align}
{\it Step 4. Passage to the limit.}
Considering estimates  (\ref{est_psik_H1}), (\ref{lipschitz_est_V2_psim}) and 
(\ref{lipschitz_est_V2_psimt}), we deduce	that the sequences $\left(\psi_k\right)_k$ and $\left(\tfrac{\partial \psi_k }{\partial t}\right)_k$ are uniformly bounded in $L^\infty(I;W)$ and $L^2(I;V)$, respectively. Then there exists a subsequence, still indexed by $k$, and function $\psi\in  H^2(\Omega)$ such that
	$$ \psi_k \longrightarrow  \psi \qquad \mbox{weakly* in} \ L^\infty(I;H^2(\Omega)),$$
	$$ \tfrac{\partial \psi_k }{\partial t}\longrightarrow  \tfrac{\partial \psi}{\partial t} \qquad \mbox{weakly in} \ L^2(I;H^1(\Omega)).$$
These convergence results imply that $\psi\in C(\bar I; H^1(\Omega))$ and as 
	$\psi_k(0)$ converges to zero in $H^1(\Omega)$, we deduce that $\psi$ verifies the initial condition
	$\psi(0)=0.$
By passing to the limit in $(\ref{faedo_galerkin_adj})_1$, we obtain for every $j\geq 1$
$$\begin{array}{ll}
\left(\tfrac{\partial \psi(t) }{\partial t}, e_j\right)+2\alpha
\displaystyle\left(D\tfrac{\partial\psi(t) }{\partial t},
De_j\right)+2\nu\left(D\psi(t),De_j\right)+
b\left(\psi(t),\boldsymbol e_j, \sigma(\tilde y(t)) \right)-b\left(\boldsymbol e_j,\psi(t), \sigma(\tilde y(t)) \right)\vspace{2mm}\\
\ \quad \qquad \qquad+
b\left(\psi(t),\tilde y(t), \sigma(\boldsymbol e_j) \right)-b\left(\tilde y(t),\psi(t), \sigma(\boldsymbol e_j) \right)
=\left(\tilde f(t), e_j\right)\end{array}$$
and by density we prove that $ \psi$ is a solution of (\ref{form_var_lin_adj_psi}).
 Moreover, $\psi$ satisfies 
$$\left\|\psi\right\|_{L^\infty(I;L^2)}^2+2\alpha
\left\|D\psi\right\|_{L^\infty(I;L^2)}^2\leq
\left\| \tilde f\right\|^2_{2,Q}
\displaystyle
e^{cT \left(1+\left(1+\alpha \right)\left\|\tilde y\right\|_{L^\infty(I;H^3)}\right)},$$		
$$2\left\|D\psi\right\|_{L^\infty(I;L^2)}^2+\alpha 
\left\|\mathbb{A}\psi\right\|_{L^\infty(I;L^2)}^2$$
$$\leq \left(\tfrac{1}{\nu}\left\| \tilde f\right\|^2_{2,Q}+cT(1+\alpha)\left\|\tilde y\right\|_{L^\infty(I;H^3)}
\left\|\psi\right\|^2_{L^\infty(I;L^2)}\right)e^{
\frac{c(1+\alpha)T}{\alpha}\left\|\tilde y\right\|_{L^\infty(I;H^3)}}.$$	
$$\left\|\tfrac{\partial \psi}{\partial t}\right\|_{2,Q}^2+
\alpha\left\|D\tfrac{\partial \psi}{\partial t}\right\|_{2,Q}^2$$
$$\leq 
c\left(\|\tilde f\|_{2,Q}^2+T\left(\nu^2+\tfrac{(1+\alpha)^2}{\alpha}\left\|\tilde y\right\|_{L^\infty(I;H^3)}^2\right)
	\left\|\psi\right\|_{L^\infty(I;H^2)}^2\right),$$
where $c$ is a constant only depending on $\Omega$. Since (\ref{adjointF}) is linear, the uniquess result is direct consequence of the first estimate.\vspace{1mm}\\
Summarizing, we have established that (\ref{adjointF}) admits a unique solution in the sense of (\ref{form_var_lin_adj_psi}) and that this solution satisfies the previous estimates. As a consequence, it follows that problem  (\ref{adjoint}) admits a unique solution in the sense of 
(\ref{form_var_lin_adj}). By taking into account the fact that
	$$\begin{array}{lll}\|\tilde f\|_{2,Q}=\left\|f\right\|_{2,Q}, & \left\|\tilde y\right\|_{L^\infty(I;H^3)}=\left\|y\right\|_{L^\infty(I;H^3)}, &\vspace{2mm}\\
	\left\|\psi\right\|_{L^\infty(I;L^2)}=\left\|p\right\|_{L^\infty(I;L^2)}, & \left\|D\psi\right\|_{L^\infty(I;L^2)}=\left\|Dp\right\|_{L^\infty(I;L^2)}, & \left\|\mathbb{A}\psi\right\|_{L^\infty(I;L^2)}=\left\|\mathbb{A}p\right\|_{L^\infty(I;L^2)}\vspace{2mm}\\
	\left\|\tfrac{\partial \psi}{\partial t}\right\|_{2,Q}=
	\left\|\tfrac{\partial p}{\partial t}\right\|_{2,Q}, & 	\left\|D\tfrac{\partial \psi}{\partial t}\right\|_{2,Q}=
	\left\|D\tfrac{\partial p}{\partial t}\right\|_{2,Q}, &
	\end{array}$$
we deduce that $p$ satisfies the claimed estimates.$\hfill\Box$\vspace{2mm}\\
We finish this section by establishing a relation between the adjoint state and the solution of a linearized state equation.
\begin{proposition} \label{green_formula}Let $y$ be in $L^\infty(I;W\cap H^3(\Omega))$ and 
$f,w\in L^2(Q)$. Then we have
	$$\int_0^T\left(w(t), p(t)\right)dt=\int_0^T\left(f(t), z(t)\right)dt,$$
where $p$ is the solution of $(\ref{adjoint})$ and where $z$ is the solution of the linearized system $(\ref{linearized})$ corresponding to $y_1=y_2=y$.
\end{proposition}
{\bf Proof.} Let $p_k(t)=\psi_k(T-t)$, where $\psi_k$ is the solution of 
(\ref{faedo_galerkin_adj}). Then $p_k=\sum_{i=1}^k \zeta_{i}^\ast(t) e_i$ satisfies
	\begin{equation}\label{pk_zk}\left\{\begin{array}{llll}
	\displaystyle-\left(\tfrac{\partial p_k(t) }{\partial t}, 
	e_j\right)-2\alpha\left(D\tfrac{\partial p_k(t) }{\partial t},
	De_j\right)+2\nu\left(Dp_k(t),De_j\right)\vspace{2mm}\\
	\ \quad \qquad \qquad+
	b\left(p_k(t),e_j, \sigma(y(t)) \right)
	-b\left(e_j,p_k(t), \sigma(y(t)) \right)\vspace{2mm}\\
	\ \quad \qquad \qquad+
	b\left(p_k(t),y(t), \sigma(e_j) \right)-b\left(y(t),p_k(t), 
	\sigma(e_j) \right)=
	\left(f(t), e_j\right) \qquad 1\leq j\leq k,\vspace{2mm}\\
	p_k(T)=0.
	\end{array}\right.\end{equation}
Let $z_k(t)=\sum_{i=1}^k \zeta_{i}^{\ast\ast}(t) e_i$ 
be the solution of (\ref{faedo_galerkin_lin}) corresponding to 
$y_1=y_2=y$, that is $z_k(t)$ satisfies
	\begin{equation}\label{zk_pk}
	\left\{\begin{array}{llll}
	\displaystyle	\left(\tfrac{\partial z_k(t) }{\partial t}, e_j\right)+2\alpha
	\displaystyle\left(D\tfrac{\partial z_k(t) }{\partial t},
	De_j\right)+2\nu\left(Dz_k(t),De_j\right)\vspace{2mm}\\
	\qquad \qquad \,\quad  
	+b\left(e_j,y(t),\sigma(z_k(t)) \right)-b\left(y(t),e_j,\sigma(z_k(t)) \right)\vspace{2mm}\\
\ \quad \qquad \qquad+
b\left(e_j,z_k(t), \sigma(y(t)) \right)-b\left(z_k(t),e_j, \sigma(y(t)) \right)=
	\left(w(t), e_j\right) \qquad 1\leq j\leq k,\vspace{2mm}\\
	z_k(0)=0.
	\end{array}\right.\end{equation}
Multiplying $(\ref{pk_zk})_1$ by $\zeta_j^{\ast\ast}(t)$
 and summing, we obtain
	\begin{align}\label{green1}
	&\displaystyle-\left(\tfrac{\partial p_k(t) }{\partial t}, 
	z_k(t)\right)-2\alpha\left(D\tfrac{\partial p_k(t) }{\partial t},
	Dz_k(t)\right)+2\nu\left(Dp_k(t),Dz_k(t)\right)\nonumber\\
	&+b\left(p_k(t),z_k(t), \sigma(y(t)) \right)
	-b\left(z_k(t),p_k(t), \sigma(y(t)) \right)\nonumber\\
	&+b\left(p_k(t),y(t), \sigma(z_k(t)) \right)-b\left(y(t),p_k(t), 
	\sigma(z_k(t)) \right)=
	\left(f(t), z_k(t)\right)\end{align}
Similarly, by multiplying $(\ref{zk_pk})_1$ by $\zeta_j^{\ast}(t)$ and summing, we obtain
	\begin{align}\label{green2}
	&\displaystyle	\left(\tfrac{\partial z_k(t) }{\partial t},p_k(t)\right)+2\alpha
	\displaystyle\left(D\tfrac{\partial z_k(t) }{\partial t},
	Dp_k(t)\right)+2\nu\left(Dz_k(t),Dp_k(t)\right)\nonumber\\
	&+b\left(p_k(t),y(t),\sigma(z_k(t)) \right)-b\left(y(t),p_k(t),\sigma(z_k(t)) \right)\nonumber\\
&+b\left(p_k(t),z_k(t), \sigma(y(t)) \right)-b\left(z_k(t),p_k(t), \sigma(y(t)) \right)=
	\left(w(t), p_k(t)\right).\end{align}
On the other hand, by applying standard rule of integration by parts we obtain
	$$\begin{array}{ll}&-\left(\tfrac{\partial p_k(t) }{\partial t}, 
	z_k(t)\right)-2\alpha\left(D\tfrac{\partial p_k(t) }{\partial t},
	Dz_k(t)\right)\vspace{2mm}\\
	&=-\left(\tfrac{\partial p_k(t) }{\partial t},
	\sigma\left(z_k(t)\right)\right)=\left(p_k(t),\tfrac{\partial \sigma\left(z_k(t)\right)}{\partial t}\right)-\tfrac{d}{dt}\left(p_k(t),\sigma\left(z_k(t)\right)\right)\vspace{2mm}\\
	&=\left(p_k(t),\tfrac{\partial \sigma\left(z_k(t)\right)}{\partial t}\right)-\tfrac{d}{dt}\left(\left(p_k(t),z_k(t)\right)+2\alpha
	\left(Dp_k(t),Dz_k(t)\right)\right)\vspace{2mm}\\
	&=\left(p_k(t),\tfrac{\partial z_k(t)}{\partial t}\right)+2\alpha\left(Dp_k(t),D\tfrac{\partial z_k(t)}{\partial t}\right)-\tfrac{d}{dt}\left(\left(p_k(t),z_k(t)\right)+2\alpha
	\left(Dp_k(t),Dz_k(t)\right)\right).
\end{array}$$
By integrating and taking into account $(\ref{pk_zk})_2$
 and $(\ref{zk_pk})_2$ we deduce that
	$$\int_0^T\left(-\left(\tfrac{\partial p_k(t) }{\partial t}, 
	z_k(t)\right)-2\alpha\left(D\tfrac{\partial p_k(t) }{\partial t},
	Dz_k(t)\right)\right)dt$$
	\begin{equation}\label{green3}=\int_0^T\left(\left(p_k(t),\tfrac{\partial z_k(t)}{\partial t}\right)+2\alpha\left(Dp_k(t),D\tfrac{\partial z_k(t)}{\partial t}\right)\right).\end{equation}
Combining (\ref{green1})-(\ref{green3}), we finally obtain
	$$\int_0^T\left(f(t), z_k(t)\right)dt=\int_0^T\left(w(t), p_k(t)\right)dt.$$
The claimed result follows by passing to the limit. $\hfill\Box$	
\section{Proof of the main results}
\setcounter{equation}{0}
{\bf Proof of the existence of an optimal control for $(P)$.}\vspace{2mm}\\
 Let $(u_{k},y_{k})_k$ be a minimizing sequence. Since $(u_{k})_k$ is uniformly bounded in the closed convex set $U_{ad}\subset L^2\left(I;H^1(\Omega)\right)$, by taking into account (\ref{state_est3})-(\ref{state_t}), we deduce that the sequence $(y_{k})_k$ is  bounded in $L^\infty\left(I; W\cap H^3(\Omega)\right))\cap H^{1}
\left(I; V\right))$. We may then extract a subsequence, still indexed by $k$, such that $(u_k)_k$ weakly converges to some $u\in U_{ad}$ in $L^2\left(I;H({\rm curl};\Omega)\right)$, $\left(y_k\right)_{k}$ weakly* converges to some $y$ in the weakly* topology of $L^\infty\left(I; H^3(\Omega)\right)\cap H^1 \left(I; V\right)$ and, by using compactness results on Sobolev spaces, in the strong topology of $C(\bar I;H^1(\Omega))$.\vspace{1mm}\\ Since $y_k(0)=y_0$, it follows that $y(0)=y_0$. Moreover, by passing to the limit in the variational formulation corresponding to $y_{k}$, we obtain
	$$\left(\tfrac{\partial y(t)}{\partial t},\phi\right)
	+	2\alpha \left(D\tfrac{\partial y(t)}{\partial t}, D\phi\right)
	+2\nu \left(Dy(t),D\phi\right)+
	b\left( \phi, y(t), \sigma(y(t))\right)-b\left(y(t),
	 \phi, \sigma(y(t))\right)=(u(t),\phi)$$
for all 
	$\phi \in V$, implying that $(u,y)$ satisfies (\ref{equation_etat}).
From the convexity and continuity of $J$, it follows the lower semicontinuity of $J$ in the weak topology and 
	$$J(u,y)\leq \liminf_kJ(u_{k},y_{k})=\inf(P),$$
showing that $(u,y)$ is a solution for $(P)$.$\hfill\Box$ \vspace{2mm}\\
{\bf Proof of the necessary optimality conditions  for $(P)$.}\vspace{2mm}\\
Taking into account Propositions \ref{existence_state},  \ref{ex_uniq_lin} and  \ref{ex_uniq_adj} and Proposition \ref{Gateau}, it follows that the corresponding state, linearized state, adjoint state exist and are unique and that the control-to-state mapping is G\^ateaux differentiable at $\bar u$. For $\rho\in ]0,1[$ and $v\in U_{ad}$, let $u_{\rho}=\bar u+\rho (v-\bar u)$ and $y_\rho$ the corresponding solution of (\ref{equation_etat}). Since $(\bar u,\bar y)$ is an optimal solution and $(u_{\rho},y_{\rho})$ is admissible, we have 
	$$\displaystyle\lim_{\rho\rightarrow 0}
	\tfrac{J(u_{\rho},y_{\rho})-
	J(\bar u,\bar y)}{\rho}\geq 0.$$
By taking into account Proposition \ref{Gateau}, we deduce that
	\begin{equation}\label{zv}\int_0^T\left(z_{\bar u,v-\bar u}(t),\bar y(t)-y_d(t)\right)+
	\lambda\left(\bar u(t),v(t)-\bar u(t)\right)dt\geq 0,\end{equation}
where $z_{\bar u,v-\bar u}$ is the (unique) solution of the linearized equation
	$$\left\{
  \begin{array}{ll}
     \displaystyle\tfrac{\partial \sigma(\boldsymbol z)}{\partial t}-\nu\Delta  \boldsymbol z+ \boldsymbol{curl}\,\sigma(\boldsymbol z)\times \bar{\boldsymbol y}+
	 \boldsymbol{curl}\,\sigma(\bar{\boldsymbol y})\times \boldsymbol z+\nabla \pi=\boldsymbol v-\bar{\boldsymbol u}&\quad\mbox{in} \ Q,\vspace{2mm}\\
	 \mathrm{div} \,  \boldsymbol z=0&\quad\mbox{in} \ Q,\vspace{2mm}\\
    \boldsymbol z\cdot \boldsymbol n=0, \qquad 
	\left(\boldsymbol n\cdot D\boldsymbol z\right)\cdot \boldsymbol\tau=0&\quad\mbox{on}\ \Sigma,\vspace{2mm}\\
	\boldsymbol z(0)=0 &\quad\mbox{in}\ \Omega.
	  \end{array}
\right .$$
Let $\bar{p}$ be the unique solution of (\ref{adj_opt_eq_alpha}). By taking into account Proposition \ref{green_formula}, we deduce that
\begin{equation}\label{cond_opt_u}\int_0^T\left(\bar{y}(t)-y_d(t),z_{\bar u,v-\bar u}(t)\right)dt=\int_0^T\left(v(t)-\bar{u}(t),\bar{p}(t)\right)dt.
	\end{equation}
The result follows by combining (\ref{zv}) and (\ref{cond_opt_u}).$\hfill\Box$


\begin{thebibliography}{99}

\bibitem{A12} {\sc N. Arada}, {\em Optimal control of shear-thinning fluids}, SIAM J. Control Optim. 50 (2012), pp. 2515-2542.

\bibitem{A13} {\sc  N. Arada}, {\em Optimal control of shear-thickening flows}, SIAM J. Control Optim. 51 (2013), 
pp. 1940-1961.

\bibitem{A14} {\sc  N. Arada}, {\em Optimal control of evolutionary quasi-Newtonian fluids}, SIAM J. Control Optim. 52 (2014), 
pp. 3401-3436.

\bibitem{BI06} {\sc A. V. Busuioc, D. Iftimie},  {\em A non-Newtonian fluid with Navier boundary conditions}, J. Dynam. Diff. Eq. 18 (2006), pp. 357-379.

\bibitem{BILN12}  {\sc A. V. Busuioc, D. Iftimie, M. C. Lopes Filho, H. J. Nussenzveig Lopes}, {\em Incompressible Euler as a limit of complex fluid models with Navier boundary conditions}, J. Differential Equations 252 (2012), pp. 624-640.


\bibitem{BR03} {\sc A. V. Busuioc, T. S. Ratiu},  {\em The second grade fluid and averaged Euler equations with Navier-slip boundary conditions},
 Nonlinearity 16 (2003), pp. 1119-1149.
 
 \bibitem{CH93}  {\sc R. Camassa, D. D. Holm}, {\em 
 An integrable shallow water equation with peaked solitons}, Phys. Rev. Lett. 71 (1993), pp. 1661-1664.
 

\bibitem{CG97}  {\sc D. Cioranescu, V. Girault}, {\em Weak and classical solutions of a family of second grade fluids}, Int. J. Nonlinear Mech. 32 (1997), pp. 317-335.

\bibitem{CO84}  {\sc D. Cioranescu, E. H. Ouazar}, {\em Existence and uniqueness for fluids of second grade}, Nonlinear Partial Differential Equations and Their Applications (Coll\`ege de France Seminar, Paris, 1982/1983), 4 (Boston, MA: Pitman) (1984), pp. 178-197.


\bibitem{CMR98}  {\sc T. Clopeau, A. Mikelic, R. Robert}, {\em On the vanishing viscosity limit for the $2D$ incompressible Navier-Stokes equations with the friction type boundary conditions}, Nonlinearity 11 (1998), pp. 1625-1636.

\bibitem{C96}  {\sc J. M. Coron}, {\em On the controllability of the $2D$ incompressible Navier-Stokes equations with the Navier-slip
boundary conditions}, ESAIM Control Optim. Calc. Var. 1 (1996), pp. 35-75.



\bibitem{FHT2}  {\sc C. Foias, D. D. Holm, E. S. Titi}, {\em The Navier-Stokes-$\alpha$ model of fluid turbulence}, Physica D 153 (2001), pp. 505-519.

\bibitem{FHT1}  {\sc C. Foias, D. D. Holm, E. S. Titi}, {\em The three dimensional viscous Camassa-Holm equations and their relation to the Navier-Stokes equations and the turbulence theory}, J. Dynam. Diff. Eq. 14 (2002), pp. 1-35.

\bibitem{GS94}  {\sc G. P. Galdi, A. Sequeira}, {\em  Further existence results for classical solutions of the equations of a second-grade fluid}, Arch. Rational Mech. Anal. 128 (1994), pp. 297-312.

\bibitem{HMR981}  {\sc D. D. Holm, J. E. Marsden, R. S. Ratiu}, {\em Euler-Poincar\'e models of ideal fluids with nonlinear dissipation}, Phys. Rev. Lett. 349 (1998), pp. 4173-4177.

 \bibitem{HMR98}  {\sc D. D. Holm, J. E. Marsden, R. S. Ratiu}, {\em The Euler-Poincar\'e equations and semi-direct products with applications to continuum theories}, Adv. Math. 137 (1998), pp. 1-81. 


\bibitem{IP06} {\sc D. Iftimie, G. Planas}, {\em Inviscid limits for the Navier-Stokes equations with Navier friction boundary conditions}, Nonlinearity 19 (2006), pp. 899-918.

\bibitem{K06}  {\sc J. P. Kelliher}, {\em Navier-Stokes equations with Navier boundary conditions for a bounded domain in the plane}, SIAM J. Math. Anal. 38 (2006), pp. 210-232.

\bibitem{LT10}  {\sc J. S. Linshiz, E. S. Titi}, {\em On the convergence rate of the Euler-$\alpha$, an inviscid second-grade complex fluid, model to Euler equations}, J. Stat. Phys. 138 (2010), pp. 305-332.

\bibitem{LNP05}  {\sc M.C. Lopes Filho, H.J. Nussenzveig Lopes, G. Planas}, {\em On the inviscid limit for $2D$ incompressible flow with Navier friction
condition}, SIAM J. Math. Anal. 36 (2005), pp. 1130-1141.

\bibitem{LNTZ15}  {\sc M.C. Lopes Filho, H. J. Nussenzveig Lopes, E. S. Titi, A. Zang}, {\em Convergence of the $2D$ Euler-$\alpha$ to Euler equations in the Dirichlet case: Indifference to boundary layers}, Physica D 292-293 (2015) pp. 51-61. 


\bibitem{O81}  {\sc E. H. Ouazar}, {\em Sur les Fluides de Second Grade}, Th\`ese 3\`eme Cycle, Universit\'e Pierre et Marie Curie, 1981.

  
\bibitem{S73} {\sc V. E. \v{S}\v{c}adilov, V. A. Solonikov}, {\em On a boundary value problem for a stationary system of Navier-Stokes
equations}, Proc. Steklov Inst. Math. 125 (1973), pp. 186-199.
























\end{thebibliography}
\end{document}